\documentclass{article}

\usepackage{proof}
\usepackage{latexsym}

\usepackage{amssymb}
\usepackage{cite}

\newcommand{\calF}{{\cal F}}
\newcommand{\calL}{{\cal L}}
\newcommand{\BS}{{\sf BS}}

\newcommand{\blem}{\begin{lemma}}
\newcommand{\elem}{\end{lemma}}
\newcommand{\bth}{\begin{theorem}}
\newcommand{\ethm}{\end{theorem}}
\newcommand{\benu}{\begin{enumerate}}
\newcommand{\eenu}{\end{enumerate}}
\newcommand{\bdes}{\begin{description}}
\newcommand{\edes}{\end{description}}
\newcommand{\bdf}{\begin{definition}}
\newcommand{\edf}{\end{definition}}
\newcommand{\bcor}{\begin{cor}}
\newcommand{\ecor}{\end{cor}}
\newcommand{\bprp}{\begin{proposition}}
\newcommand{\eprp}{\end{proposition}}

\newtheorem{theorem}{Theorem}[section]
\newtheorem{definition}[theorem]{Definition}
\newtheorem{proposition}[theorem]{Proposition}
\newtheorem{lemma}[theorem]{Lemma}
\newtheorem{cor}[theorem]{Corollary}

\newcommand{\bprf}{{\bf Proof}.\hspace{2mm}}

\newcommand{\eprf}{\hspace*{\fill} $\Box$}

\newcommand{\beqn}{\begin{equation}}
\newcommand{\eeqn}{\end{equation}}
\newcommand{\beqnarr}{\begin{eqnarray}}
\newcommand{\eeqnarr}{\end{eqnarray}}
\newcommand{\beqnarrs}{\begin{eqnarray*}}
\newcommand{\eeqnarrs}{\end{eqnarray*}}

\newcommand{\alp}{\alpha}

\newcommand{\del}{\delta}
\newcommand{\Del}{\Delta}
\newcommand{\ome}{\omega}
\newcommand{\Ome}{\Omega}
\newcommand{\bet}{\beta}
\newcommand{\gam}{\gamma}
\newcommand{\Gam}{\Gamma}

\newcommand{\Sig}{\Sigma}
\newcommand{\tht}{\theta}

\newcommand{\Lam}{\Lambda}
\newcommand{\vphi}{\varphi}

\newcommand{\fal}{\forall}
\newcommand{\exi}{\exists}

\newcommand{\Rarw }{\Rightarrow}

\newcommand{\lrarw}{\leftrightarrow}
\newcommand{\Lrarw}{\Leftrightarrow}

\newcommand{\la}{\langle}
\newcommand{\ra}{\rangle}
\newcommand{\lc}{\lceil}
\newcommand{\rc}{\rceil}
\newcommand{\setm}{\setminus}
\newcommand{\spand}{\,\&\,}
\newcommand{\restrict}{\!\upharpoonright\!}

\begin{document}

\title{Intuitionistic fixed point theories over set theories}


\author{Toshiyasu Arai
\\
Graduate School of Science,
Chiba University
\\
1-33, Yayoi-cho, Inage-ku,
Chiba, 263-8522, JAPAN
\\
tosarai@faculty.chiba-u.jp
}

\maketitle

\begin{abstract}
In this paper we show that the intuitionistic fixed point theory $\mbox{FiX}^{i}(T)$ over
set theories $T$ is a conservative extension of $T$
if $T$ can manipulate finite sequences and has the full foundation schema.
\end{abstract}

\section{Intuitionistic fixed point theory over set theories $T$}

For a theory $T$ in a laguage $\calL$, let $\mathcal{Q}(X,x)$ be an $X$-positive formula 
in the language $\calL\cup\{X\}$ with an extra unary predicate symbol $X$.
Introduce a fresh unary predicate symbol $Q$ together with the axiom stating that $Q$ is a fixed point of $\mathcal{Q}(X,x)$:
\beqn\label{eq:fixax}
\fal x[Q(x)\lrarw \mathcal{Q}(Q,x)]
\eeqn
By the completeness theorem, it is obvious that the resulting extension of $T$ is conservative over $T$,
though it has a non-elementary speed-up over $T$
when $T$ is a recursive theory containing the elementary recursive arithmetic {\sf EA},
cf. \cite{speedup}.

When $T$ has an axiom schema, e.g., $T={\sf PA}$, the Peano arithmetic with the complete induction schema,
let us define the fixed point extension $\mbox{FiX}({\sf PA})$ to have the induction schema for any formula with the fixed point predicate $Q$.
Then $\mbox{FiX}({\sf PA})$ is stronger than {\sf PA}, e.g., $\mbox{FiX}({\sf PA})$ proves the consistency of {\sf PA}.
For the proof-theretic strength of the fixed point theory $\mbox{FiX}({\sf PA})$, see \cite{Feferman, attic, Avigad}.

On the other side, 
W. Buchholz\cite{BuchholzAML97} shows that an {\it intuitionistic\/} fixed point theory 
over the intuitionistic (Heyting)  arithmetic {\sf HA} for strongly positive formulae $\mathcal{Q}(X,x)$
 is proof-theoretically reducible to {\sf HA}.
In a language of arithmetic strongly positive formulae with respect to $X$ 
are generated from arithmetic formulae and atomic ones $X(t)$ by means of positive connectives $\lor,\land,\exi,\fal$.
Then R\"uede and Strahm\cite{Strahm} extends the result to the intuitionistic fixed point theory 
$\mbox{FiX}^{i}({\sf HA})$
for strictly positive formulae $\mathcal{Q}(X,x)$, 
in which the predicate symbol $X$ does not occur in the antecedent $\vphi$ of implications $\vphi\to\psi$ 
nor in the scope of negations $\lnot$.
Indeed
as shown in \cite{intfix} 
$\mbox{FiX}^{i}({\sf HA})$ is
a conservative extension of {\sf HA}.

However this might mislead us.
Namely one might think that the conservation holds for the fixed point extensions because
the theory $T={\sf HA}$ is intuitionistic.
Actually this is not the case.
For example, the intuitionistic fixed point theory $\mbox{FiX}^{i}({\sf PA})$
over the {\it classical\/} arithmetic {\sf PA} is a conservative extension of {\sf PA}.
For, if $\mbox{FiX}^{i}({\sf PA})$ proves an arithmetical sentence $A$, then
$\mbox{FiX}^{i}({\sf HA})$ proves $B\to A$ for a {\sf PA}-provable sentence $B$.
Since $\mbox{FiX}^{i}({\sf HA})$ is conservative over {\sf HA},
we see that {\sf HA} proves $B\to A$, and ${\sf PA}\vdash A$.

Our proof in \cite{intfix} is a proof-theoretic one by showing that the fixed point axiom (\ref{eq:fixax}) is
eliminable quickly.
The crux is that the {\it underlying logic\/} is intuitionistic.
\\

\noindent
{\bf Digression}.
Let $\widehat{ID}^{i}(\mbox{acc})$ be an intuitionistic theory obtained from $\widehat{ID}^{i}(\mbox{strict})=\mbox{FiX}^{i}({\sf HA})$
by restricting $\mathcal{Q}(X,x)$ to {\it accessible formulas\/}, i.e.,
$\mathcal{Q}(X,x)\equiv(A(x)\land\fal y(B(y,x)\to X(y)))$ for arithmetical formulas $A,B$.
The following Lemma \ref{lem:Strahm} is shown in \cite{Strahm}.
\blem\label{lem:Strahm}{\rm (\cite{Strahm})}
\benu
\item\label{lem:Strahm1}
$\widehat{ID}^{i}(\mbox{{\rm strict}})$ is conservative over $\widehat{ID}^{i}(\mbox{{\rm acc}})$ with respect to almost negative formulas.

\item\label{lem:Strahm2}
The classical theory $\widehat{ID}(\mbox{{\rm acc}})$ is interpretable in the classical arithmetic {\sf PA}.
\eenu
\elem
Lemma \ref{lem:Strahm}.\ref{lem:Strahm1} is shown by a recursive realizability interpretation following Buchholz\cite{BFPS},
and the interpretation in the proof of Lemma \ref{lem:Strahm}.\ref{lem:Strahm2} is done by a diagonalization argument.
Specifically it is observed that there is an arithmetical fixed point for accessible operators, classically.
Then they conclude that $\widehat{ID}^{i}(\mbox{strict})$ is conservative over the intuitionistic arithmetic {\sf HA}
with respect to negative formulas.

Let us try to prove the full conservation result in \cite{intfix} along the line in \cite{Strahm}.
The intuitionistic version of Lemma \ref{lem:Strahm}.\ref{lem:Strahm2} is easy to see, which says that
$\widehat{ID}^{i}(\mbox{{\rm acc}})$ is a conservative extension of {\sf HA}.
 Let $A(x)$ and $B(y,x)$ be arithmetical formulae.
Let $y<_{B}x:\Lrarw B(y,x)$ and $y\leq_{B}^{*}x$ denote its reflexive and transitive closure.
Namely
$y\leq_{B}^{*}x$ iff there exists a non-empty sequence $(x_{n},\ldots,x_{0})$ such that
$x_{n}=y$, $x_{0}=x$ and $\fal i<n(x_{i+1}<_{B}x_{i})$.
Then $A^{*}(x):\Lrarw\fal y\leq_{B}^{*}x\, A(y)$ is an arithmetical fixed point for accessible operators 
$\mathcal{Q}(X,x)\equiv(A(x)\land\fal y(B(y,x)\to X(y)))$ provably in {\sf HA}, i.e.,
${\sf HA}\vdash\fal x[A^{*}(x)\lrarw(A(x)\land\fal y<_{B}x\, A^{*}(y))]$.
The problem is to extend Lemma \ref{lem:Strahm}.\ref{lem:Strahm1} to all arithmetical formulae,
which means that $\widehat{ID}^{i}(\mbox{{\rm strict}})$ is conservative over  $\widehat{ID}^{i}(\mbox{{\rm acc}})$
with respect to any arithmetical formulae.
If a combination of realizability interpretation and forcing works as in \cite{Beeson1979}, then it would yield 
the full conservativity.
However it is hard to show the soundness of the forcing stating that
if $\widehat{ID}^{i}(\mbox{{\rm acc}}){\bf a}\vdash A$, then 
$\widehat{ID}^{i}(\mbox{{\rm acc}})\vdash \fal p\exi q\supset p(q\Vdash A)$
since $p\Vdash \fal y(B(y,x)\to Q(y))$ is not an accessible formula, but strictly positive.
\\

\noindent
In this paper we extend the observation in \cite{intfix}  for set theories $T$.

Let $T$ be a set theory in the language $\{\in,=\}$.

Fix an $X$-strictly positive formula $\mathcal{Q}(X,x)$ in the language $\{\in,=,X\}$ with an extra unary predicate symbol $X$.
In $\mathcal{Q}(X,x)$ the predicate symbol $X$ occurs only strictly positive.
The language of $\mbox{FiX}^{i}(T)$ is $\{\in,=,Q\}$ with a fresh unary predicate symbol $Q$.
The axioms in $\mbox{FiX}^{i}(T)$ consist of the following:
\benu
\item
All provable sentences in $T$ (in the language $\{\in,=\}$).

\item
Foundation schema for any formula $\vphi$ in the language $\{\in,=,Q\}$:
\beqn\label{eq:Qind}
\fal x(\fal y\in x\,\vphi(y)\to\vphi(x))\to\fal x\,\vphi(x)
\eeqn

\item
Fixed point axiom (\ref{eq:fixax}).
\eenu

The underlying logic in $\mbox{FiX}^{i}(T)$ is defined to be the intuitionistic (first-order predicate) logic (with equality).
$\fal x,y(x=y\to Q(x)\to Q(y))$ is an axiom.

In this paper we show the following Theorem \ref{th:consvintfix} for a weak base set theory $\BS$ defined in the next section \ref{sect:BS}.

\bth\label{th:consvintfix}
$\mbox{{\rm FiX}}^{i}(T)$ is a conservative extension of any set theory $T\supset\BS$.
\end{theorem}

We need Theorem \ref{th:consvintfix} in \cite{liftupK, liftupZF, cnsvrfijslfix, liftupKn}
for proof-theoretic analyses of set theories for weakly compact cardinals, first-order reflecting ordinals,  {\sf ZF} and
second-order indescribable cardinals.
In these analyses, a provability relation $\mathcal{H}\vdash^{\alp}_{c}\Gam$ derived from operator controlled derivations
is defined to be a fixed point of a strictly positive formula.

Let us mention the contents of the paper.
In section \ref{sect:BS} a weak base theory $\BS$ is introduced, and it is shown that $\BS$ can manipulate finite sequences
and partially define truth.
In section \ref{sect:code} a class of codes $Code$ and a binary relation $\prec$ on it are defined,
and it is shown that the transfinite induction schema with respect to $\prec$ is provable in $\BS$ up to {\it each\/} code.
The order type of the well founded relation $\prec$ is the next epsilon number to the order type of the class of ordinals
in the universe.
In section \ref{sect:finan} a sequent calculus for $\mbox{{\rm FiX}}^{i}(T)$ is introduced, and in section \ref{sect:prflemma}
Theorem \ref{th:consvintfix} is proved by a finitary analysis of the proofs in the sequent calculus for $\mbox{{\rm FiX}}^{i}(T)$.

\section{Basic set theory \BS}\label{sect:BS}

In this section we introduce a basic set theory $\BS$, and show that $\BS$ can manipulate finite sequences of sets,
 thereby can encode syntax, and define truth partially.
 
Consider the following functions $\calF_{i}\,(i<9)$,
$\calF_{0}(x,y)=\{x,y\}$,
$\calF_{1}(x,y)=\cup x$,
$\calF_{2}(x,y)=x\setm y$,
$\calF_{3}(x,y)=\{u\cup \{v\} : u\in x, v\in y\}$,
$\calF_{4}(x,y)=dom(x)=\{u\in\cup\cup x : \exi v\in\cup\cup x(\la u,v\ra\in x)\}$,
$\calF_{5}(x,y)=rng(x)=\{v\in\cup\cup x : \exi u\in\cup\cup x(\la u,v\ra\in x)\}$,
$\calF_{6}(x,y)=\{\la v,u\ra\in y\times x: v\in u\}$,
$\calF_{7}(x,y)=\{\la u,v,w\ra : \la u,v\ra\in x, w\in y\}$, and
$\calF_{8}(x,y)=\{\la u,w,v\ra : \la u,v\ra\in x, w\in y\}$, where
$\la v,u\ra=\{v,\{v,u\}\}$ and $\la u,v,w\ra=\la u,\la v,w\ra\ra$.

Note that each $\calF_{i}$ is {\it simple\/} in the sense that for any $\Del_{0}$-formula $\vphi(z)$,
$\vphi(\calF_{i}(x,y))$ is $\Del_{0}$.
For each $i$, $\calF_{i}(x,y,z)$ denotes a $\Del_{0}$-formula stating $\calF_{i}(x,y)=z$.

\bdf\label{df:BS}
{\sf BS} {\rm is the set theory in the language} $\{\in,=\}$.
{\rm Its axioms are  Extensionality,  
Foundation schema, and}
$\{\fal x,y\exi z\, \calF_{i}(x,y,z): i<9\}$.

{\rm A set-theoretic function} $f:V^{n}\to V$ {\rm is} $\Sig_{1}^{\BS}${\it -definable\/}
{\rm if there exists a} $\Sig_{1}${\rm -formula} $\vphi(x_{1},\ldots,x_{n},y)$ {\rm for which}
$\BS\vdash\fal x_{1},\ldots,x_{n}\exi ! y\, \vphi(x_{1},\ldots,x_{n},y)$, {\rm and}
$f(x_{1},\ldots,x_{n})=y$ {\rm iff} $V\models\vphi(x_{1},\ldots,x_{n},y)$.

{\rm A relation} $R\subset V^{n}$ {\rm is} $\Del_{1}^{\BS}$ {\rm if there exist}
$\Sig_{1}${\rm -formulae} $\vphi,\psi$ {\rm such that}
\\
$\BS\vdash\fal x_{1},\ldots,x_{n}[\vphi(x_{1},\ldots,x_{n})\lrarw\lnot\psi(x_{1},\ldots,x_{n})]$,
{\rm and}
$(x_{1},\ldots,x_{n})\in R$ {\rm iff} $V\models\vphi(x_{1},\ldots,x_{n})$.
\edf

A formula is said to be $\Del_{0}(\ome)$ iff every quantifier occurring in it is
a bounded quantifier $\exi m<n,\fal m<n$ with bounds $n\in\ome$.

\bprp\label{prp:BS}
\benu
\item\label{prp:BS1}
The Cartesian product $a\times b$ is a $\Sig_{1}^{\BS}$-function with a $\Del_{0}$-graph.

\item\label{prp:BS2}
$\BS$ proves $\Del_{0}$-Separation:
$\BS\vdash\fal a,b\exi c[c=\{x\in a:\vphi(x,b)\}]$ for each $\Del_{0}$-formula $\vphi$.

\item\label{prp:BS3}
$\ome\times V\ni(n,a)\mapsto{}^{<n}a,{}^{n}a$ are $\Sig_{1}^{\BS}$-functions.
$z={}^{<n}a$ is $\Del_{0}$, and $z={}^{n}a$ is $\Del_{1}^{\BS}$.

\item\label{prp:BS4}
The class of $\Del_{1}^{\BS}$-relations is closed under propositional connectives $\lnot,\lor$
and bounded quantifications $\exi m<n,\fal m<n$ with bounds $n\in\ome$.

For $\Sig_{1}^{\BS}$-functions $f$ and $\Del_{0}(\ome)$-formula $\tht(y)$,
$\tht(f(\vec{x}))$ is $\Del_{1}^{\BS}$.

\item\label{prp:BS5}
The class of $\Sig_{1}^{\BS}$-functions is closed under compositions and primitive recursion on $\ome$.
The latter means that if $g: V^{n}\to V$ and $h:\ome\times V^{n+1}\to V$
 are $\Sig_{1}^{\BS}$-functions, then so is the function $f:\ome\times V^{n}\to V$ defined by
$f(0,\vec{x})=g(\vec{x})$ and $f(n+1,\vec{x})=h(n,\vec{x},f(n,\vec{x}))$ for $\vec{x}=x_{1},\ldots,x_{n}$.

\item\label{prp:BS6}
For the transitive closure $trcl(a)$, $x\in trcl(a)$ is $\Del_{1}^{\BS}$.
\eenu
\eprp
\bprf
\ref{prp:BS}.\ref{prp:BS1}.
Let $G(a)=\calF_{3}(\{\emptyset\},a)=\{\{x\}: x\in a\}$.
Then by $\calF_{3}(G(a),b)=\{\{x,y\} : x\in a, y\in b\}$,
we have
$a\times b=\{\{\{x\},\{x,y\}\} : x\in a, y\in b\}=\calF_{3}(G(G(a)),\calF_{3}(G(a),b))$.
$a\times b=c$ is a $\Del_{0}$-formula.
\\

\noindent
\ref{prp:BS}.\ref{prp:BS2}. 
Standard, cf. \cite{barwise}, pp. 63-67 using Proposition \ref{prp:BS}.\ref{prp:BS1}.
\\

\noindent
\ref{prp:BS}.\ref{prp:BS3}.
Noting
${}^{n+1}a=\{x\cup \{y\}  : x\in {}^{n}a, y\in \{n\}\times a\}=\calF_{3}({}^{n}a, \{n\}\times a)$,
$\BS$ proves the existence of ${}^{n}a$ by induction on $n\in\ome$.
Next observe that ${}^{<n}a=z$ is $\Del_{0}$ since $x\in{}^{<n}a$ as well as $x\in{}^{m}a$ is $\Del_{0}$ and for $n>0$,
${}^{<n}a= z$ iff $z\subset{}^{<n}a$ and 
$\{\emptyset\}={}^{0}a\subset z$ and $\fal m<n-1\fal x\in z\cap{}^{m}a\fal b\in a[x\cup\{\la m,b\ra\}\in z]$.
Therefore $z={}^{n}a$ iff $z=(({}^{<n+1}a)\setm({}^{<n}a))$.
\\

\noindent
\ref{prp:BS}.\ref{prp:BS4}.
For $\Sig_{1}$-formula $\exi x\,\tht(m,x)$ with $\Del_{0}$-matrix $\tht$,
$\BS$ proves that $\fal m<n\exi x\,\tht(m,x)\lrarw\exi y\fal m<n\exi x\in y\,\tht(m,x)$
by induction on $n\in\ome$.
\\

\noindent
\ref{prp:BS}.\ref{prp:BS5}.
Let the function $f$ be defined from $\Sig_{1}^{\BS}$-functions $g,h$ by 
$f(0,\vec{x})=g(\vec{x})$ and $f(n+1,\vec{x})=h(n,\vec{x},f(n,\vec{x}))$.
Then
$f(n,\vec{x})=y$ iff there exists a function $F$ with $dom(F)=n+1$ such that
$F(0)=g(\vec{x})$, $\fal i<n[F(i+1)=h(i,\vec{x},F(i))]$ and $y=F(n)$.
By induction on $n\in\ome$ $\BS$ proves $\fal n\in\ome\fal\vec{x}\exi! y[f(n,\vec{x})=y]$.
Moreover $f(n,\vec{x})=y$ is $\Sig_{1}$ by Proposition \ref{prp:BS}.\ref{prp:BS4}.
\\

\noindent
\ref{prp:BS}.\ref{prp:BS6}.
From Proposition \ref{prp:BS}.\ref{prp:BS5} we see that
$(n,a)\mapsto\cup^{(n)}a$ is a $\Sig_{1}^{\BS}$-function,
where $\cup^{(0)}a=a$ and $\cup^{(n+1)}a=\cup(\cup^{(n)}a)$.
Hence $x\in trcl(a)\Lrarw \exi n\in\ome(x\in \cup^{(n)}a)\Lrarw \fal b(\cup b\subset b \land a\subset b \to x\in b)$.
\eprf
\\

\noindent
From Proposition \ref{prp:BS} we see that $\BS$ can encode syntax, e.g., formulae in the language $\{\in,=\}$.
Let $\lc Fml\rc\subset\ome$ denote the set of codes $\lc\vphi\rc$ of formulae $\vphi$ in $\{\in,=\}$.

We can assume that $\lc Fml\rc$ is $\Del_{1}^{\BS}$,
and manipulations on it, e.g.,
$(\lc\vphi\rc,\lc\psi\rc)\mapsto\lc\vphi\lor\psi\rc$,
$\lc\vphi\lor\psi\rc\mapsto\la\lc\vphi\rc,\lc\psi\rc\ra$,
are all $\Sig_{1}^{\BS}$.
Moreover for $x\in \lc Fml\rc$,
let $var(x)$ denote the set $\{n\in\ome: v_{n} \mbox{ occurs freely in } x\}$, and
$ass(x,y)$ the set of function $f:var(x)\to y$.
Both $x\mapsto var(x)$ and $(x,y)\mapsto ass(x,y)$ are $\Sig_{1}^{\BS}$-functions.
Let $\models\lc\vphi\rc[a]$ denote the satisfaction relation for formulae $\vphi$
and $a\in ass(\lc\vphi\rc,y)$ for a $y$.

For formula $\vphi$ in $\{\in,=\}$, $\lc Sbfml\rc(\vphi)$ denotes the finite set of codes of subformulae of $\vphi$.

\blem\label{lem:truthdef}
For {\rm each} formula $\vphi$ in the language $\{\in,=\}$, the satisfaction relation
$\{(x,a): x\in \lc Sbfml\rc(\vphi), a\in ass(x), \models x[a]\}$ for subformulae of $\vphi$ is 
$\BS$-definable in such a way that $\BS$ proves that
 $\vphi(v_{0},\ldots,v_{m-1})\lrarw \models\lc\vphi(v_{0},\ldots,v_{m-1})\rc[a]$
 for $a(i)=v_{i}$,
 $\models\lc\vphi_{0}\lor\vphi_{1}\rc[a]\lrarw \models\lc\vphi_{0}\rc[a_{0}]\lor \models\lc\vphi_{1}\rc[a_{1}]$
 for $a_{i}=a\restrict var(\lc\vphi_{i}\rc)$ and subformulae $\vphi_{i}$,
 $\models\lc\exi v_{m}\,\vphi\rc[a]\lrarw \exi b[\models\lc\vphi\rc[a\cup\{\la m,b\ra\}]$ for subformula $\exi v_{m}\,\vphi$,
 and similarly for $\land,\fal$.
\elem
\bprf
It suffices to $\Del_{1}^{\BS}$-define the satisfaction relation for subformulae of a given $\Del_{0}$-formula $\vphi$.
This is seen as in \cite{SchindlerZeman}, p.613
using Propositions \ref{prp:BS}.\ref{prp:BS3} and \ref{prp:BS}.\ref{prp:BS4}.
Note that we don't need the existence of transitive closures to bound range $y$ of the assignments $a:var(x)\to y$
since there are only finitely many subformulae of the given $\vphi$.
\eprf

\section{Codes}\label{sect:code}

Let us define a class $Code$ of codes and a binary relation $\prec$ on it recursively.
It is shown that the transfinite induction schema with respect to $\prec$ is provable in $\BS$ up to {\it each\/} code.

The class $Code$ of codes together with the relation $\prec$ is essentially a notation system of `ordinals'
whose order type is the next epsilon number to the order type of the class of ordinals
in the universe $V$.
To define such a notation system, we need at least ordinal addition $\alp+\bet$ 
and exponentiation with base, say $\ome$, $\ome^{\alp}$ at hand.
However $\BS$ is too weak to $\Del_{1}$-define $\alp+\bet$ and $\ome^{\alp}$,
since it lacks  $\Del_{0}$-Collection.
In other words, the order type $\Lam$ of the class of ordinals in the well founded universe $V\models\BS$ need not to be an epsilon number
nor even an additive principal number, which is closed under $\alp+\bet$.
Indeed, $L_{\alp}\models\BS$ for any limit ordinal $\alp$.

Instead of $\Del_{0}$-Collection,
we collect formal expressions called {\it products\/} $\bar{a}_{1}\times\cdots\times\bar{a}_{n}$ of codes $\bar{a}_{i}$ for 
$a_{i}\in V\cup\{V\}$ first,
and then collect formal expressions called {\it sums\/} $\alp_{1}\#\cdots\#\alp_{n}$ of products $\alp_{i}$.
Intuitively $\#$ denotes the natural (commutative) sum, and $\times$ the natural product,
if the code $\bar{a}$ is replaced by the ordinal $2 rank(a)$.
Each sum is defined to be smaller than a code $\Ome$, 
which is interpreted as the least additive principal number $(\Lam+1)^{\ome}$ above $\Lam$.
Then introduce formal expressions $\Ome^{\alp}\bet$, which is intended to be an exponential function.
These three operations $\times,\#$ and $(\alp,\bet)\mapsto\Ome^{\alp}\bet$ on codes are needed in the ordinal assignment
to proofs defined in Definition \ref{df:ordass}.
The relation $\prec$ on codes is well founded, but not a linear ordering.
For our proof-theoretic analysis, the linearity of $\prec$ is dispensable, 
the base $\Ome$ can be replaced by $2$, and
$\bar{a}$ by $rank(a)$.
Definitions \ref{df:sump} and \ref{df:codeless} simplify the matters.

First let us define a class $Sum$ and a relation $\prec_{p}$ on it.
$\ell(\alp)$ is the {\it length\/} of $\alp\in Sum$.

\bdf\label{df:sum}
{\rm Let} $\bar{a}:=\la 0,a\ra$ {\rm for} $a\in V$, {\rm and} $\bar{V}:=\la 1,0\ra$.
$\ell(\bar{a}):=0$.
\benu
\item
{\rm A} {\it product\/} {\rm is either} $\bar{1}$ {\rm or}
$\bar{a}_{1}\times\cdots\times\bar{a}_{n}:=\la 2, \bar{a}_{1},\ldots,\bar{a}_{n}\ra$ {\rm for}
$a_{1},\ldots,a_{n}\in V\cup\{V\}$ {\rm with} $a_{i}\neq 0,1$ {\rm and} $n>0$.
$Prod$ {\rm denotes the class of all products.}

$\ell(\bar{a}_{1}\times\cdots\times\bar{a}_{n})=\max\{\ell(\bar{a}_{1}),\ldots,\ell(\bar{a}_{n})\}+1$.
{\rm When} $n=0$, {\rm let}
$\alp_{1}\times\cdots\times\alp_{n}:=\bar{1}$.

\item
{\rm A} {\it sum of products\/} {\rm is a set}
$\alp_{1}\#\cdots\#\alp_{n}:=\la 3, \alp_{1},\ldots,\alp_{n}\ra$ {\rm with} $\alp_{i}\in Prod$ {\rm and} $n\geq 0$.
$Sum$ {\rm denotes the class of all sums of products.}

$\ell(\alp_{1}\#\cdots\#\alp_{n})=\max\{\ell(\alp_{1}),\ldots,\ell(\alp_{n})\}+1$.

{\rm When} $n=0$, {\rm let}
$\alp_{1}\#\cdots\#\alp_{n}:=\bar{0}$.
\eenu
$Prod$ {\rm is a subclass of} $Sum$.
\edf

Let us introduce some operations and `computation rules' on sums.
\benu

\item

$\times$ and $\#$ are defined to be commutative, i.e.,
$\alp_{1}\#\cdots\#\alp_{n}=\alp_{\pi(1)}\#\cdots\#\alp_{\pi(n)}$ and
$ \alp_{1}\times\cdots\times \alp_{n}= \alp_{\pi(1)}\times\cdots\times \alp_{\pi(n)}$
for any permutation $\pi\in n!$.

These means that $\alp_{1}\#\cdots\#\alp_{n}$ and $\alp_{1}\times\cdots\times \alp_{n}$
are actually multisets of products and codes $\bar{a}$.

\item
$\bar{0}$ is the zero element.
$\alp_{1}\#\cdots\#\alp_{n}\#\bar{0}=\alp_{1}\#\cdots\#\alp_{n}$,
and
$\alp_{1}\times\cdots\times \alp_{n}\times\bar{0}:=\bar{0}$.

\item

$\bar{1}$ is the unit for $\times$,
$ \alp_{1}\times\cdots\times \alp_{n}\times\bar{1}:= \alp_{1}\times\cdots\times \alp_{n}$.

\item
Associative laws
$(\alp_{1}\#\cdots\#\alp_{n})\#(\bet_{1}\#\cdots\#\bet_{m}):=\alp_{1}\#\cdots\#\alp_{n}\#\bet_{1}\#\cdots\#\bet_{m}$
for $\{\alp_{1},\ldots,\alp_{n}\}\cup\{\bet_{1},\ldots,\bet_{m}\}\subset Prod$,
and $(\alp_{1}\times\cdots\times \alp_{n})\times\bet:=\alp_{1}\times\cdots\times \alp_{n}\times\bet$
for $\{\alp_{1},\ldots,\alp_{n}\}\cup\{\bet\}\subset \{\bar{a}: a\in V\}\cup\{\bar{V}\}$.

\item
Distributive laws
$(\alp_{1}\#\cdots\#\alp_{n})\times\bet:=(\alp_{1}\times\bet)\#\cdots\#(\alp_{n}\times\bet)$
where either $\{\alp_{1},\ldots,\alp_{n}\}\subset Prod$, $\bet\in \{\bar{a}: a\in V\}\cup\{\bar{V}\}$.
\eenu
Therefore any combination of products by $\#$ and $\times$ is equal to (reduced to) a sum of products.

\bdf\label{df:sump}
$\alp\prec_{p}\bet$ {\rm for} $\alp,\bet\in Sum$ {\rm is defined recursively as follows.}
{\rm Let}
$\alp\preceq_{p}\bet:\Lrarw (\alp\prec_{p}\bet)\lor(\alp=\bet)$.
\benu
\item\label{df:sump0}
$\bar{0}\prec_{p} \alp$ {\rm for any sum} $\alp\neq\bar{0}$.
\item\label{df:sump1}
{\rm For} $\gam\in Prod$, 
$\{\gam_{i}: 0\leq i\leq n\}\subset Prod\cup\{\bar{0}\}$ {\rm and} $n>0$,
\[
a\in^{(n)} b \in V\cup\{V\} \spand \fal i\leq n(\gam_{i}\preceq_{p}\gam) \Rarw (\gam_{0}\times\bar{a})\#\gam_{1}\#\cdots\#\gam_{n}\prec_{p}\gam\times\bar{b}
\]
{\rm where}
$a\in^{(n+1)}b:\Lrarw a\in \cup^{(n)}b$ {\rm with} $\cup^{(0)}b=b$ {\rm and} $\cup^{(n+1)}b=\cup(\cup^{(n)}b))$.

\item\label{df:sump2}
$\alp_{0}\prec_{p}\alp_{1} \spand \bet_{0}\preceq_{p}\bet_{1} \Rarw \alp_{0}\#\bet_{0}\prec_{p}\alp_{1}\#\bet_{1}$.
\eenu
\edf
Definition \ref{df:sump}.\ref{df:sump1} says that
if $a=a_{n}\in a_{n-1}\in \cdots\in a_{1}\in b$, then
$(\gam\times a)\#(\gam\cdot n)\prec_{p}\gam\times b$ for $\gam\cdot n=\gam\#\cdots\#\gam\, (n \mbox{ times }\gam)$.

\bprp
The relation $\prec_{p}$ on $Sum$ is transitive.
\eprp

Next let us define the class of codes $Code$.
$\ell(\alp)$ is the {\it length\/} of $\alp\in Code$.

\bdf\label{df:code}
 \benu

 \item
 $\bar{0}\in Code$.
 $\ell(\bar{0})=0$.
 
 \item
 $PCode\subset Code$.
 
 \item
 $\alp\in Code \spand \bar{0}\neq\bet\in Sum \Rarw 
\Ome^{\alp}\bet:=\la 4,\alp,\bet\ra\in PCode$.

$\ell(\Ome^{\alp}\bet)=\max\{\ell(\alp),\ell(\bet)\}+1$.

 \item
 $\alp_{1},\ldots,\alp_{n}\in PCode \spand n>0 \Rarw\alp_{1}\#\cdots\#\alp_{n}:=\la 3,\alp_{1},\ldots,\alp_{n}\ra\in Code$.

 $\ell(\alp_{1}\#\cdots\#\alp_{n})=\max\{\ell(\alp_{1}),\ldots,\ell(\alp_{n})\}+1$.
 \eenu

\edf

Let us introduce some operations and `computation rules' on codes.
\benu

\item
Again $\#$ is defined to be commutative, and
$\bar{0}$ is the zero element.

\item

$\bar{1}$ is the unit,
 $\Ome^{\bet}:=\Ome^{\bet}\bar{1}$,
$\Ome^{\bar{0}}=\bar{1}$ and $\alp_{1}\times\cdots\times \alp_{n}:=\alp_{1}\times\cdots\times \alp_{n}\times\bar{1}$
for $\alp_{i}\in \{\bar{a}: a\in V\}\cup\{\bar{V}\}$.
Also $\alp=\Ome^{\bar{0}}\alp$.
Thus $Sum\subset Code$.
\item
{\rm When} $n=1$, $\alp_{1}\#\cdots\#\alp_{n}$  {\rm is identified with} $\alp_{1}\in PCode$.

\item
Exponential law
$\Ome^{\gam}(\Ome^{\bet}\alp):=\Ome^{\gam\#\bet}\alp$ for $\alp\in Sum$.

\item
Associative laws for $\#$ and
Distributive laws
$\Ome^{\bet}(\alp_{1}\#\cdots\#\alp_{n}):=\Ome^{\bet}\alp_{1}\#\cdots\#\Ome^{\bet}\alp_{n}$
where
$\{\alp_{1},\ldots,\alp_{n}\}\subset PCode$.

\eenu

Next we define a binary relation $\prec$ on $Code$ recursively as follows.

\bdf\label{df:codeless}
\benu
\item
$\bar{0}\prec \alp$ {\rm for any code} $\alp\neq\bar{0}$.

\item
{\rm Let} $\bar{0}\not\in\{\bet_{i}:i<n\}\cup\{\bet\}\subset Sum$ {\rm and} 
$\{\alp_{i}:i<n\}\cup\{\alp\}\subset Code$ {\rm with} $\alp_{i}\neq\alp_{j}\,(i\neq j)$.
\[
\fal i<n[(\alp_{i},\bet_{i})\prec_{lex}(\alp,\bet)] \Rarw \Ome^{\alp_{0}}\bet_{0}\#\cdots\#\Ome^{\alp_{n-1}}\bet_{n-1}\prec\Ome^{\alp}\bet
\]
{\rm where for codes} $\alp,\gam\in Code$ {\rm and sums} $\bet,\del\in Sum$
\[
(\alp,\bet)\prec_{lex}(\gam,\del):\Lrarw \alp\prec\gam \mbox{ {\rm or }} (\alp=\gam\spand \bet\prec_{p}\del)
.\]

\item
$\alp_{0}\prec\alp_{1} \spand \bet_{0}\preceq\bet_{1} \Rarw \alp_{0}\#\bet_{0}\prec\alp_{1}\#\bet_{1}$.
\eenu

\edf


The following Proposition \ref{prp:code} is easily seen.

\bprp\label{prp:code}
\benu
\item\label{prp:code0}
Both $Code$ and $\prec$ are $\Del_{1}^{\BS}$.

\item\label{prp:code1}
The relation $\prec$ on $Code$ is transitive.

\item\label{prp:code2}
For $\alp,\bet\in Sum$, $\alp\prec_{p}\bet\Lrarw\alp\prec\bet$,
where $\alp=\Ome^{\bar{0}}\alp$.

\item\label{prp:code4}
$\bar{0}$ is the least element.

\item\label{prp:code5}
For $a\in b \in V\cup\{V\}$, sums of products $\gam,\gam^{\prime},\del$ 
$\gam^{\prime}\preceq\gam  \Rarw (\gam\times\bar{a})\#\gam^{\prime}\prec(\gam\times\bar{b})\#\del$.

\item\label{prp:code6}
 $\bet\prec \alp\#\bet$ if $\alp\neq \bar{0}$.
 
\item\label{prp:code7}
$\gam\prec\alp\#\bet \Rarw \gam\prec\alp \lor \exi\bet_{0}\prec\bet(\gam=\alp\#\bet_{0})$.

\item\label{prp:code75}
$\del\prec\bar{a}_{1}\times\cdots\times\bar{a}_{n}\Rarw \exi i\leq n\exi b\in a_{i}[\del\preceq (\bar{b}\prod_{j\neq i}\bar{a}_{j})\#\prod_{j\neq i}\bar{a}_{j}]$
for $\del\in Sum$ and $a_{i}\in V\cup\{V\}$.

\item\label{prp:code8}
Both $(\alp,\bet)\mapsto\alp\#\bet$ and
$(\alp,\bet)\mapsto \Ome^{\alp}\bet$ are monotonic in each argument.

\item\label{prp:code9}
$\alp,\bet\prec\Ome^{\alp}\bet$ if $\alp,\bet\neq \bar{0}$.

\item\label{prp:code10}
$\alp_{1},\alp_{2}\prec\bet \Rarw \Ome^{\alp_{1}}\#\Ome^{\alp_{2}}\prec \Ome^{\bet}$.

\item\label{prp:code11}
$
\bet_{0}\prec\bet 
\Rarw
\Ome^{\bet_{0}}(\alp\#\alp)\prec\Ome^{\bet}\alp$.
\eenu
\eprp

For a binary relation $<$ and formulae $\vphi$, let
\beqnarrs
Prg[\vphi,<] & :\Lrarw & \fal x(\fal y<x\,\vphi(y)\to \vphi(x))
\\
TI[\vphi,<,a] & :\Lrarw & Prg[\vphi,<] \to \fal x<a\,\vphi(x)
\\
TI[<,a] & := & \{TI[\vphi,<,a] : \vphi \mbox{ is a formula}\}
\eeqnarrs

$T\vdash TI[<,a]$ means that $T\vdash TI[\vphi,<,a]$ for any formula $\vphi$, and
$T\vdash TI[<,a]\to TI[<,b]$ means that
for any $\vphi$ there exists a formula $\psi$ such that
$T\vdash TI[\psi,<,a]\to TI[\vphi,<,b]$.

\blem\label{lem:wfveps}
For {\rm each} code $\alp\in Code$,
$
\BS\vdash TI[\prec,\alp]
$.
\elem

Lemma \ref{lem:wfveps} is shown by metainduction on the length $\ell(\alp)$ of codes $\alp$
using the following Proposition \ref{prp:wfveps}.

\bprp\label{prp:wfveps}
\benu
\item\label{prp:wfveps1}
$\BS\vdash TI[\prec_{p},\alp]\land TI[\prec_{p},\bet]\to TI[\prec_{p},\alp\#\bet]$.
Similarly for $\prec$.

\item\label{prp:wfveps2}
$\BS\vdash TI[\prec_{p},\bar{V}]$, i.e.,
$\BS\vdash Prg[\vphi,\prec_{p}]\to\fal a\in V\cup\{V\}\vphi(\bar{a})$ for any formula $\vphi$.

\item\label{prp:wfveps25}
For any formula $\vphi$,
$\BS\vdash \fal n<\ome\fal \{a_{i}\}_{i<n}\subset V\cup\{V\}
[\fal i<n (\fal x\prec\prod_{j\neq i}\bar{a}_{j}\,\vphi(x)\land \fal y\prec_{p}\bar{a}_{i}\fal x\prec y\prod_{j\neq i}\bar{a}_{j}\,\vphi(x) )
\to \fal x\prec\prod_{i<n}\bar{a}_{i}\,\vphi(x)]$.

\item\label{prp:wfveps3}
$\BS\vdash \fal\alp\in Sum\,TI[\prec_{p},\alp]$, i.e,
$\BS\vdash Prg[\vphi,\prec_{p}]\to \fal \alp\in Sum\,\vphi(\alp)$ for any formula $\vphi$.

\item\label{prp:wfveps45}
$\BS\vdash Prg[\vphi,\prec]\land \fal \alp_{0}\prec\alp\fal \gam\fal x\prec\Ome^{\alp_{0}}\gam\, \vphi(x) \land
 \fal \gam\prec_{p}\bet\fal x\prec\Ome^{\alp}\gam\,\vphi(x) \to \fal x\prec\Ome^{\alp}\bet\,\vphi(x)$.

\item\label{prp:wfveps4}
$\BS\vdash TI[\prec,\alp]  \to TI[\prec,\Ome^{\alp}\bet]$.
\eenu
\eprp
\bprf
\ref{prp:wfveps}.\ref{prp:wfveps1}.
This follows from the fact $\BS\vdash Prg[\vphi,\prec_{p}]\to Prg[\vphi_{\alp\#},\prec_{p}]$
for $\vphi_{\alp\#}(x):\Lrarw \vphi(\alp\#x)$ using Proposition \ref{prp:code}.\ref{prp:code7}.
\\

\noindent
\ref{prp:wfveps}.\ref{prp:wfveps2}.
This is seen from Foundation schema.
\\

\noindent
\ref{prp:wfveps}.\ref{prp:wfveps25}.
This is seen from Propositions \ref{prp:code}.\ref{prp:code75} and \ref{prp:wfveps}.\ref{prp:wfveps1}.
\\

\noindent
\ref{prp:wfveps}.\ref{prp:wfveps3}.
By Proposition \ref{prp:wfveps}.\ref{prp:wfveps1} it suffices to show 
$\BS\vdash Prg[\vphi,\prec_{p}]\to \fal \alp\in Prod\,\vphi(\alp)$ for any formula $\vphi$.
We show this by induction on the number $n$ of components in products $\bar{a}_{1}\times\cdots\times\bar{a}_{n}$.
The case $n=1$, $Prg[\vphi,\prec_{p}]\to \fal a\in V\cup\{V\}\, \vphi(\bar{a})$ follows from Proposition \ref{prp:wfveps}.\ref{prp:wfveps2}.
Let $Prod_{n}$ denote the class of all products such that the number of components is at most $n$.
Suppose $Prg[\vphi,\prec_{p}]$ and $ \fal \alp\in Prod_{n}\,\vphi(\alp)$ for a formula $\vphi$.
Let $a_{i}\in V\cup\{V\}$.
Then by Proposition \ref{prp:wfveps}.\ref{prp:wfveps25} we have 
$\fal i<n+1\fal y\prec_{p}\bar{a}_{i}\fal x\prec_{p} y\prod_{j\neq i}\bar{a}_{j}\,\vphi(x)  \to \fal x\prec_{p}\prod_{i<n+1}\bar{a}_{i}\,\vphi(x)$.
In other words,
$\fal i<n\fal y\prec_{p}\bar{a}_{i}\fal x\prec_{p} y\prod_{j\neq i}\bar{a}_{j}\,\vphi(x)  \to Prg[\vphi_{n},\prec_{p}]$,
where
$\vphi_{n-k}(y):\Lrarw \fal x\prec_{p}y\prod_{j\neq n-k}\bar{a}_{j}\,\vphi(x)$.
Thus by \ref{prp:wfveps}.\ref{prp:wfveps2} we have
$\fal i<n\fal y\prec_{p}\bar{a}_{i}\fal x\prec_{p} y\prod_{j\neq i}\bar{a}_{j}\,\vphi(x)  \to \fal x\prec_{p}\prod_{i<n+1}\bar{a}_{i}\,\vphi(x)$.
In this way we see
$\fal i<n+1-k\fal y\prec_{p}\bar{a}_{i}\fal x\prec_{p} y\prod_{j\neq i}\bar{a}_{j}\,\vphi(x)  \to \fal x\prec_{p}\prod_{i<n+1}\bar{a}_{i}\,\vphi(x)$
by induction on $k\leq n+1$.
Hence $\fal x\prec_{p}\prod_{i<n+1}\bar{a}_{i}\,\vphi(x)$, i.e., $ \fal \alp\in Prod_{n+1}\,\vphi(\alp)$.
\\

\noindent
\ref{prp:wfveps}.\ref{prp:wfveps45}.
This is seen from Proposition \ref{prp:code}.\ref{prp:code75} and Definition \ref{df:codeless}.
\\

\noindent
\ref{prp:wfveps}.\ref{prp:wfveps4}.
Suppose $Prg[\vphi,\prec]$ and
$\fal \alp_{0}\prec\alp\fal \bet\fal x\prec\Ome^{\alp_{0}}\bet\, \vphi(x)$.
Then by Proposition \ref{prp:wfveps}.\ref{prp:wfveps45}
we have $Prg[\vphi_{\Ome^{\alp}},\prec_{p}]$, where
$\vphi_{\Ome^{\alp}}(\bet):\Lrarw \fal x\prec\Ome^{\alp}\bet\, \vphi(x)$.
Hence by Proposition \ref{prp:wfveps}.\ref{prp:wfveps3}
$\fal \bet\fal x\prec\Ome^{\alp}\bet\, \vphi(x)$.
Thus we have shown
$Prg[\vphi,\prec]\to Prg[\mbox{j}[\vphi],\prec]$, where
$\mbox{j}[\vphi](\alp_{0}):\Lrarw \fal \bet\fal x\prec\Ome^{\alp_{0}}\bet\, \vphi(x)$.
Hence by $TI[\prec,\alp]$ we have
$\fal \bet\fal x\prec\Ome^{\alp}\bet\, \vphi(x)$.
\eprf
\\

Lemma \ref{lem:wfveps} is now seen by metainduction on the length $\ell(\alp)$ of codes $\alp$
using Propositions \ref{prp:wfveps}.\ref{prp:wfveps1}, \ref{prp:wfveps}.\ref{prp:wfveps3} and \ref{prp:wfveps}.\ref{prp:wfveps4}.

\section{Finitary analysis of $\mbox{{\rm FiX}}^{i}(T)$}\label{sect:finan}

When the set theory $T$ is sufficiently strong, e.g., when $T$ comprises Kripke-Platek set theory,
we could prove Theorem \ref{th:consvintfix} as in \cite{intfix}, i.e., 
first the finitary derivations of set-theoretic sentences $\vphi$ in $\mbox{{\rm FiX}}^{i}(T)$ 
are embedded to infinitary derivations of a sequent
$\tht\Rarw\vphi$ for a provable sentence $\tht$ in $T$, then partial cut-elimination is possible.
This results in a $\Del_{1}$-definable infinitary derivation of the same sequent $\tht\Rarw\vphi$
 in which there
occur no fixed point formulae.
The depth of the derivation is bounded by an exponential ordinal tower.
Then transfinite induction shows that $\tht\Rarw\vphi$ is true.
By formalizing the infinitary arguments straightforwardly  in $T$ 
we would see that the end formula $\vphi$ is true in $T$.
To formalize the infinitary analysis in a weaker theory $T$, 
we need a finitary treatment of it as in \cite{BuchholzAML91}.

Let us take another route in terms of Gentzen-Takeuti's finitary analyses 
of finite derivations as in \cite{PT2} since its formalization in a weak (set) theory is a trivial matter.

In what follows we work in a set theory $T\supset\BS$.

$\alp,\bet,\gam,\ldots$ range over codes in $Code$,
while $a,b,c,\ldots$ over sets in the universe $V$.
$A,B,C,\ldots$ denote formulae in the language 
$\calL_{V}:=\{\in,=,Q\}\cup\{\bar{a}:a\in V\}$,
where $\bar{a}:=\la 0,a\ra$ is the name (individual constant) for the set $a$.
A term is either a name or a variable.
$\iota,\nu,\ldots$ denote terms.

Let us introduce a sequent calculus for 
transfinite induction schema (\ref{eq:Qind}) and the fixed point axiom (\ref{eq:fixax}).
Logical connectives are $\lor,\land,\to, \exists,\forall$.
$\lnot A:\equiv(A\to\bot)$.

A {\it sequent\/} is a pair of a finite set $\Gam$ of formulae, and a formula $A$, denoted
$\Gamma\Rightarrow A$.
Its intended meaning is the implication $\bigwedge\Gam\to A$.
$\Gam$ is the {\it antecedent\/}, and $A$ the {\it succedent\/} of the sequent $\Gam\Rarw A$.
For finite sets $\Gam$, $\Del$ and a formula $A$,
$\Gam,\Del:=\Gam\cup\Del$ and $\Gam,A:=\Gam\cup\{A\}$.

$\bot$ stands ambiguously for false atomic sentences $\bar{a}\in\bar{b}$ for $a\not\in b$,
and $\bar{a}=\bar{b}$ for $a\neq b$.

The {\it initial sequents\/} are
\[
\Gam, \iota=\nu,A(\iota)\Rarw A(\nu)
\, ; \mbox{\hspace{5mm}} \Gamma,\bot\Rightarrow A
\]

The {\it inference rules\/} are $(LQ)$, $(RQ)$, $(L\lor)$, $(R\lor)$, $(L\land)$, $(R\land)$, $(L\to)$, $(R\to)$, 
$(L\exists)$, $(R\exists)$, $(L\forall)$, $(R\forall)$, $(cut)$, $(chain)$, $(ind)$, $(Rep)$ and $(E)$.

\[
\infer[(LQ)]
{\Gamma,Q(\iota) \Rightarrow C}
{\Gamma,Q(\iota), \mathcal{Q}(Q,\iota) \Rightarrow C}
\: ;\:
\infer[(RQ)]
{\Gamma\Rightarrow Q(\iota)}
{\Gamma\Rightarrow \mathcal{Q}(Q,\iota)}
\]

\[
\infer[(L\lor)]
{\Gamma,A_{0}\lor A_{1}\Rightarrow C}
{
\Gamma,A_{0}\lor A_{1},A_{0}\Rightarrow C
&
\Gamma,A_{0}\lor A_{1},A_{1}\Rightarrow C
}
\: ;\:
\infer[(R\lor)]
{\Gamma\Rightarrow A_{0}\lor A_{1}}
{\Gamma\Rightarrow A_{i}}
\,(i=0,1)
\]

\[
\infer[(L\land)]
{\Gamma,A_{0}\land A_{1}\Rightarrow C}
{
\Gamma,A_{0}\land A_{1},A_{i}\Rightarrow C
}
\, (i=0,1)
\: ;\:
\infer[(R\land)]
{\Gamma\Rightarrow A_{0}\land A_{1}}
{
\Gamma\Rightarrow A_{0}
&
\Gamma\Rightarrow A_{1}
}
\]

\[
\infer[(L\to)]
{\Gamma,A\to B\Rightarrow C}
{
\Gamma,A\to B\Rightarrow A
&
\Gamma,A\to B,B\Rightarrow C
}
\: ;\:
\infer[(R\to)]
{\Gamma\Rightarrow A\to B}
{\Gamma,A\Rightarrow B}
\]

\[
\infer[(L\exists)]
{\Gamma,\exists x B(x)\Rightarrow C}
{
\Gamma,\exists x B(x), B(y)\Rightarrow C
}
\: ;\:
\infer[(R\exists)]
{\Gamma\Rightarrow \exists x B(x)}
{\Gamma\Rightarrow B(\iota)}
\]
The {\it eigenvariable\/} $y$ in $(L\exi)$ does not occur in the lower sequent $\Gamma,\exists x B(x)\Rightarrow C$.

 \[
 \infer[(L\forall)]
 {\Gamma,\forall x B(x)\Rightarrow C}
 {\Gamma,\forall x B(x),B(\iota)\Rightarrow C}
 \: ;\:
 \infer[(R\forall)]
 {\Gamma\Rightarrow \forall x B(x)}
 {
 \Gamma\Rightarrow B(y)
 }
 \]
 The {\it eigenvariable\/} $y$ in $(R\fal)$ does not occur in the lower sequent $\Gamma\Rightarrow \forall x B(x)$.

  \[
  \infer[(cut)]
  {\Gamma,\Del\Rightarrow C}
  {
  \Gamma\Rightarrow A
  &
 \Del,A\Rightarrow C
  }
  \]
where $A$ is the  {\it cut formula\/} of the $(cut)$.

\[
\infer[(chain)]{\Gam,\Del\Rarw C}
{
\Gam_{k}\Rarw A_{k}
&
\cdots
&
\Gam_{1}\Rarw A_{1}
&
\Del,A_{k},\ldots, A_{1}\Rarw C
}
\]
where $\Gam=\Gam_{k}\cup\cdots\cup\Gam_{1}$, and $A_{k},\ldots, A_{1}\,(k> 0)$
is a non-empty list of {\it strictly positive\/} formulae.

The inference rule $(chain)$ is a series of several $(cut)$'s with the strictly positive {\it cut formulae\/} $A_{k},\ldots, A_{1}$.
Writing $\mbox{\boldmath$\Gam$}$ for the list $\Gam_{k},\ldots, \Gam_{1}$ and 
$\mbox{\boldmath$A$}$ for the list $A_{k},\ldots, A_{1}$, the inference rule is denoted
\[
\infer[(chain)]{\mbox{\boldmath$\Gam$},\Del\Rarw C}
{
\mbox{\boldmath$\Gam$}\Rarw \mbox{\boldmath$A$}
&
\Del,\mbox{\boldmath$A$}\Rarw C
}
\]

\[
\infer[(ind)]
{\Gam\Rarw C}
{
\Gam,\fal y\in x\, A(y)\Rarw A(x)
&
A(\iota),\Gam\Rarw C
&
\Gam\Rarw \iota\in \nu
}
\]
The {\it eigenvariable\/} $x$ does not occur in the lower sequent $\Gam\Rarw C$.

\[
\infer[(Rep)]{\Gam,\Del\Rarw A}{\Gam\Rarw A}
\: ;\:
\infer[(E)]{\Gam\Rarw A}{\Gam\Rarw A}
\]
This inference rule $(E)$ is called the height rule in \cite{pntind}, and its meaning is explained in Definition \ref{df:height}
as in \cite{BuchholzAML91}.

A {\it proof\/} in this sequent calculus is a finite labelled tree according to the above initial sequents and inference rules.
$s,t,u,\ldots$ denote the nodes in proof trees.
$s:\Gam\Rarw A$ indicates that the sequent $\Gam\Rarw A$ is the label of the node $s$.
The label $\Gam\Rarw A$ of $s$ is denoted $Seq(s)$.

Suppose that a $\{\in,=\}$-sentence $\vphi$ is provable in $\mbox{FiX}^{i}(T)$.
Then there exists a $T$-provable sentence $\tht$ 
 such that the sequent $\tht\Rarw\vphi$ is provable in the sequent calculus.
In what follows fix $\vphi,\tht$ and a proof $P_{0}$ of $\tht\Rarw\vphi$.

\bdf\label{df:purevar}
{\rm A proof in the sequent calculus is said to enjoy the} {\it pure variable condition\/}
{\rm if}
 \benu
 \item\label{df:purevar1}
  {\rm any eigenvariables (of $(L\exi), (R\fal), (ind)$) are distinct from each other,}
 \item\label{df:purevar2}
  {\rm any eigenvariable does not occur in its end sequent}, {\rm and}
 \item\label{df:purevar3}
  {\rm if a free variable occurs in an upper sequent of an inference rule but not in the lower sequent, then the variable is one  of the eigenvariables of the inference rule.} 
 \eenu
\edf

Without loss of generality we can assume that any proof enjoys the pure variable condition.
Otherwise 
rename the eigenvariables to satisfy (\ref{df:purevar1}) and (\ref{df:purevar2}) in Definition \ref{df:purevar}, then
replace the redundant free variables by an individual constant, e.g., the empty set $\bar{\emptyset}$
to satisfy (\ref{df:purevar3}).

\bdf
{\rm The} {\it end-piece\/} {\rm of a proof tree} $P$ {\rm is a collection of nodes in} $P$
{\rm such that any inference rule below it is one of} $(cut), (chain), (Rep)$ {\rm and} $(E)$.
\edf

If a proof enjoys the pure variable condition and its end sequent 
consists solely of sentences, no free variable occurs in its end-piece.

\begin{definition}\label{df:depth}
{\rm The} {\it depth} $dp(A)<\ome$ {\rm of a formula} $A$ {\rm is defined as follows.}
\benu
\item
$dp(A)=0$ {\rm if $A$ is}
$Q$-free, {\rm i.e., the fixed point predicate} Q {\rm does not occur in} $A$.

{\rm In what follows consider the case when} $Q$ {\rm occurs in} $A$.

\item
$dp(A)=2$ {\rm if} $A$ {\rm is strictly positive (with respect to $Q$)}.

{\rm In what follows consider the case when} $Q$ {\rm occurs in} $A$, {\rm and}
$A$ {\rm is not strictly positive.}

\item
$dp(A)=\max\{dp(A_{0}), dp(A_{1})\}+1$ {\rm if} $A\equiv(A_{0}\lor A_{1}), (A_{0}\land A_{1}), (A_{0}\to A_{1})$.
\item
$dp(A)=dp(A_{0})+1$ {\rm if} $A\equiv(\exi x\, A_{0}), (\fal x\, A_{0})$.
\eenu
\end{definition}
Note that $dp(A)\neq 1$.

Let $P$ be a proof in the sequent calculus, and $s$ a node in the proof tree $P$.
We assign the {\it height\/} $h(s;P)<\ome$ recursively as follows.

\bdf\label{df:height}
\benu
\item
$h(s;P)=0$ {\rm if} $Seq(s)$ {\rm is the end sequent of} $P$.

{\rm In what follows let} $Seq(s)$ {\rm be an upper sequent of an inference rule} $J$ {\rm with the lower sequent} $Seq(s_{0})$.
\item
$h(s;P)=h(s_{0};P)+1$ {\rm if} $J$ {\rm is an} $(E)$.

\item
$h(s;P)=\max\{h(s_{0};P),2\}$
{\rm if} $J$ {\rm is a} $(chain)$ {\rm with its} {\it rightmost\/} {\rm upper sequent} $Seq(s)$.
\[
\infer[(chain)]{s_{0}: \mbox{\boldmath$\Gam$},\Del\Rarw C}
{
\mbox{\boldmath$s$}: \mbox{\boldmath$\Gam$} \Rarw \mbox{\boldmath$A$}
&
s: \Del,\mbox{\boldmath$A$}\Rarw C
}
\]
\item
$h(s;P)=h(s_{0};P)$ {\rm in all other cases.}

{\rm Note that for upper sequents $\mbox{\boldmath$s$}=s_{k},\ldots, s_{1}$ of a $(chain)$ other than the rightmost one $s$, 
we have
$h(s_{i};P)=h(s_{0};P)$, i.e., the height is the same.}
\eenu
{\rm A proof} $P$ {\rm is said to be} {\it height-normal\/}
{\rm if the following four conditions hold.}
\benu
\item
{\rm For any $(chain)$ 
occurring in $P$}
\[
\infer[(chain)]{s: \mbox{\boldmath$\Gam$},\Del\Rarw C}
{
 \mbox{\boldmath$\Gam$} \Rarw \mbox{\boldmath$A$}
&
\Del,\mbox{\boldmath$A$}\Rarw C
}
\]
$h(s;P)=0$, {\rm in other words there is neither $(E)$ nor no rightmost upper sequent of $(chain)$ below any $(chain)$.}

\item
{\rm For any $(cut)$ occurring in $P$}
\[
\infer[(cut)]{s: \Gam,\Del\Rarw C}
{
\Gam\Rarw A
&
\Del,A\Rarw C
}
\]
$h(s;P)\geq dp(A)$.

\item
{\rm For any $(ind)$ occurring in $P$}
\[
\infer[(ind)]
{s: \Gam\Rarw C}
{
\Gam,\fal y\in x\, A(y)\Rarw A(x)
&
A(\iota),\Gam\Rarw C
&
\Gam\Rarw \iota\in \nu
}
\]
$h(s;P)\geq dp(\fal y\in \nu\, A(y))$.

\item
{\rm Any $(chain)$ and $(E)$ in $P$ is in the end-piece.}

\eenu

\edf

Without loss of generality we can assume that the given sequent calculus proof $P_{0}$ of $\tht\Rarw\vphi$
does not contain any $(chain)$, and
is height-normal.
Otherwise add some inference rules $(E)$ at the end of the proof.
\\

Let $P$ be a height-normal proof in the sequent calculus, and $s$ be a node in the proof tree $P$.
We assign a code $o(s;P)\in Code$ recursively as follows.
For $n>0$, $\bar{1}\cdot n:=\bar{1}\#\cdots\#\bar{1}$ with $n$ times $\bar{1}$'s.

\bdf\label{df:ordass}
\benu
\item\label{df:ordass1}
$o(s;P)=\bar{1}\cdot 2$ {\rm if} $S$ {\rm is an initial sequent.}

{\rm In what follows let} $Seq(s)$ {\rm be the lower sequent of an inference rule} $J$ {\rm with its upper sequents}
$\{s_{i}: Seq(s_{i})\}_{i<m}$.

\item\label{df:ordass2}
$o(s;P)=o(s_{0};P)\#\bar{1}$ {\rm if} $J$ {\rm is one of the inference rules} 
$(LQ)$, $(RQ)$, $(R\lor)$, $(L\land)$,  $(R\to)$, 
$(L\exists)$, $(R\exists)$, $(L\forall)$, {\rm and} $(R\forall)$.

\item\label{df:ordass3}
$o(s;P)=o(s_{0};P)\# o(s_{1};P)$
{\rm if} $J$ {\rm is one of the inference rules} $(L\lor)$, $(R\land)$, {\rm and} $(L\to)$.

\item\label{df:ordass4}
$o(s;P)=o(s_{0};P)\#o(s_{1};P)$ {\rm if} $J$ {\rm is a} $(cut)$.

\item\label{df:ordass5}
$o(s;P)=\Ome_{2}(o(s_{m-1};P))(o(s_{0};P)\#\cdots\#o(s_{m-2};P))$ {\rm if} $J$ {\rm is a} $(chain)$ {\rm where}
$\Ome_{2}(\alp):=\Ome^{\Ome^{\alp}}${\rm :}
\[
\infer[(chain)\, J]{s:\Gam,\Del\Rarw C}
{
s_{0}: \Gam_{0}\Rarw A_{0}
&
\cdots
&
s_{m-2}:\Gam_{m-2}\Rarw A_{m-2}
&
s_{m-1}:\Del,A_{0},\ldots,A_{m-2}\Rarw C
}
\]
{\rm with} $\Gam=\bigcup_{i<m-1}\Gam_{i}$.

\item\label{df:ordass6}
$o(s;P)=((o(s_{0};P)\#\bar{1}\cdot 6)\times mj(\nu))\#o(s_{1};P)\# o(s_{2};P)$
{\rm if} $J$ {\rm is an $(ind)$:}
\[
\infer[(ind)\, J]
{s:\Gam\Rarw C}
{
s_{0}: \Gam,\fal y\in x\, A(y)\Rarw A(x)
&
s_{1}: A(\iota),\Gam\Rarw C
&
s_{2}: \Gam\Rarw \iota\in\nu
}
\]
{\rm where for terms} $\nu$, 
\[
mj(\nu):=
\left\{
\begin{array}{ll}
\bar{a} & \mbox{{\rm if }} \nu=\bar{a} \mbox{ {\rm with }} a\in V
\\
\bar{V} &  \mbox{{\rm if }} \nu \mbox{ {\rm is a variable}}
\end{array}
\right.
\]

\item\label{df:ordass7}
$o(s;P)=o(s_{0};P)$ {\rm if} $J$ {\rm is a} $(Rep)$.

\item\label{df:ordass8}
$o(s;P)=\Ome^{o(s_{0};P)}$ {\rm if} $J$ {\rm is an} $(E)$.
\eenu
{\rm Finally let}
$o(P)=o(s_{end};P)$ {\rm for the end sequent} $s_{end}$ {\rm of} $P$.
\edf
The role of operations $\#,\times$ and $\Ome^{\alp}\bet$ in `ordinal' assignment $o(s;P)$ are as follows.
The sum $\alp\#\bet$ collects two subproofs together, and $\times$ is needed to multiply $\nu$ in transfinite induction $(ind)$ up to $\nu$,
cf. {\bf Case 2} in section \ref{sect:prflemma}.
Exponentiation is used first in the rule $(E)$, i.e., to measure an increase of ordinal depths in lowering cut rank,
and second in the rule $(chain)$.
The assignment $\Ome^{\Ome^{\alp}}(\alp_{k}\#\cdots\#\alp_{1})$ in $(chain)$ comes from Lemma 9 in \cite{intfix},
which in turn is inspired by the quick cut-elimination strategy in \cite{MintsFest, Mintsmono} along Kleene-Brouwer ordering of
infinitary derivations.
Lexicographic comparing, i.e., multiplication of $\Ome^{\Ome^{\alp}}$ and $\alp_{k}\#\cdots\#\alp_{1}$
is used in {\bf Case 9}, and
a doubly exponential $\Ome^{\Ome^{\alp}}$ is needed in {\bf Case 6} and {\bf Case 7},
once multiplications are introduced.
Note that when exponent $\alp$ decreases, one can duplicate multiplier $\bet$ in $\Ome^{\alp}\bet$, 
cf. Proposition \ref{prp:code}.\ref{prp:code11}.

Since any $(chain)$ and $(E)$ in $P$ is in the end-piece, $o(s;P)$ is in $Sum$
if $s$ is above the end-piece.

A formula in $\calL_{V}$ is said to be an {\it instance\/} of a formula $A$
if it is obtained from $A$ by substituting terms for free variables.

\bdf
$ISbfml(P_{0})$ {\rm denotes the class of all instances of subformulae of formulae occurring in} $P_{0}$.

{\rm Call a proof} {\it restricted\/} {\rm (with respect to $P_{0}$) if it is height-normal, enjoys the pure variable condition,
any formula occurring in it is in $ISbfml(P_{0})$, and its end sequent consists solely of 
$Q$-free sentences.}
\edf

For $\alp\in Code$ let $\tau(\alp)$ denote the formula stating that
for any restricted proof $P$
if $o(P)\preceq\alp$,
then its end sequent is true.
Note here that the satisfaction relation for the 
$Q$-free formulae in $Sbfml(P_{0})$ (the set of subformulae of formulae occurring in $P_{0}$)
 or equivalently
the partial truth definition for the 
$Q$-free sentences in $ISbfml(P_{0})$ is
{\sf BS}-definable, a fortiori $T$-definable by Lemma \ref{lem:truthdef}.

We show the following Lemma \ref{lem:main}.

\blem\label{lem:main}
$T$ proves that $\tau(\alp)$ is progressive, i.e.,
\[
T\vdash\fal\alp\in Code[\fal\bet\prec\alp \tau(\bet)\to \tau(\alp)]
.\]
\elem

Then Theorem \ref{th:consvintfix} is seen as follows.
Lemmata \ref{lem:wfveps} and \ref{lem:main} yields $\tau(o(P_{0}))$,
and hence the end sequent $\tht\Rarw\vphi$ of $P_{0}$ is true in $T$.
Therefore $T\vdash\vphi$.

\section{Proof of Lemma \ref{lem:main}}\label{sect:prflemma}

In this section we show the Lemma \ref{lem:main}.
We work in $T$.

Let $P$ be a restricted proof of a sequent $\Gam_{0}\Rarw A_{0}$.
Suppose as the IH(=Induction Hypothesis) that the end sequents of restricted proofs with smaller codes
are true.
We need to show that $\Gam_{0}\Rarw A_{0}$ is true.
It suffices to show that there are restricted proofs $P_{i}\, (i\in I)$ of sequents $S_{i}$ such that
$o(P_{i})\prec o(P)$ for any $i\in I$ and
if all of  $S_{i}$ are true, then so is $\Gam_{0}\Rarw A_{0}$.
\\

\noindent
{\bf Case 1}.
The case when there exists an initial sequent in the end-piece of $P$.

Since there are no free variables in the end-piece, any initial sequent in it is either $\Lam,\bot\Rarw A$ or
$\Lam,A\Rarw A$.

If the end sequent $\Gam_{0}\Rarw A_{0}$ itself is an initial sequent, i.e., $\{\bot,A_{0}\}\cap\Gam_{0}\neq\emptyset$,
then there is nothing to prove.
In what follows assume that this is not the case.

Consider first the case that an initial sequent $\Lam,\bot\Rarw A$ is in the end-piece.
Then the formula $\bot$ in the antecedent has to vanish somewhere as a cut formula.
Let $P$ be the following:

\[
\infer*{\Gam_{0}\Rarw A_{0}}
{
 \infer[(chain)]{s: \mbox{\boldmath$\Gam$},\Gam,\Del\Rarw C}
 {
\mbox{\boldmath$\Gam$}\Rarw \mbox{\boldmath$A$}
 &
 \infer*[Q]{s_{0}:\Gam\Rarw \bot}{}
 &
 \infer*{\Del,\mbox{\boldmath$A$},\bot\Rarw C}{\Lam,\bot\Rarw A}
 }
}
\]
Let $Q_{C}$ denote the proof obtained from the subproof $Q$ of $s_{0}:\Gam\Rarw \bot$
by replacing $\bot$ by $C$ in the succedents of
sequents $\Gam^{\prime}\Rarw\bot$ in $Q$.
Let $P^{\prime}$ be the following:
\[
\infer*{\Gam_{0}\Rarw A_{0}}
{
 \infer[(Rep)]{s: \mbox{\boldmath$\Gam$}, \Gam,\Del\Rarw C}
 {
 \infer*[Q_{C}]{s_{0}:\Gam\Rarw C}{}
 }
}
\]
Then it is clear that $P^{\prime}$ is restricted.
Moreover $o(s;P^{\prime})=o(s_{0};P^{\prime})=o(s_{0};P)\prec o(s;P)$ 
by Propositions \ref{prp:code}.\ref{prp:code6} and  \ref{prp:code}.\ref{prp:code9}.
Hence $o(P^{\prime})\prec o(P)$ by Proposition  \ref{prp:code}.\ref{prp:code8}.
From IH we see that $\Gam_{0}\Rarw A_{0}$ is true.

The case when $\bot$ vanishes at a $(cut)$ is similar.

Next consider the case that an initial sequent $\Lam,A\Rarw A$ is in the end-piece.
Then one of the formulae $A$ has to vanish somewhere as a cut formula of $J$, which is either a $(chain)$ or a $(cut)$.
Suppose $J$ is a $(chain)$, and
let $P$ be one of the followings:
\[
\infer*{\Gam_{0}\Rarw A_{0}}
{
 \infer[J]{s: \mbox{\boldmath$\Gam$},\Gam,\Del\Rarw A}
 {
 \mbox{\boldmath$\Gam$}\Rarw \mbox{\boldmath$A$}
 &
 \infer*[Q]{s_{0}:\Gam\Rarw A}{}
 &
 \infer*{s_{1}:\Del,\mbox{\boldmath$A$},A\Rarw A}{\Lam,A\Rarw A}
 }
}
\: ;\:
\infer*{\Gam_{0}\Rarw A_{0}}
{
 \infer[J]{s: \mbox{\boldmath$\Gam$},\Gam,A,\Del\Rarw C}
 {
\mbox{\boldmath$\Gam$}\Rarw \mbox{\boldmath$A$}
 &
 \infer*{s_{1}:\Gam,A\Rarw A}{\Lam,A\Rarw A}
 &
 \infer*[Q]{s_{0}:\Del,\mbox{\boldmath$A$},A\Rarw C}{}
 }
}
\]
Let $P^{\prime}$ be the followings:
\[
\infer*{\Gam_{0}\Rarw A_{0}}
{
 \infer[(Rep)]{s: \mbox{\boldmath$\Gam$},\Gam,\Del\Rarw A}
 {
 \infer*[Q]{s_{0}:\Gam\Rarw A}{}
 }
}
\: ;\:
\infer*{\Gam_{0}\Rarw A_{0}}
{
 \infer[J]{s: \mbox{\boldmath$\Gam$},\Gam,A,\Del\Rarw C}
 {
 \mbox{\boldmath$\Gam$}\Rarw \mbox{\boldmath$A$}
 &
 \infer*[Q]{s_{0}:\Del,\mbox{\boldmath$A$},A\Rarw C}{}
 }
}
\]
In the right hand side $J$ denotes two consecutive $(E)$'s if $ \mbox{\boldmath$A$}$ is the empty list, and
an $(chain)$ otherwise.
In each case $P^{\prime}$ is restricted.
Moreover $o(s_{0};P^{\prime})=o(s_{0};P)$ and $o(s_{1};P)\neq \bar{0},\bar{1}$.
Hence $o(s;P^{\prime})\prec o(s;P)$ by Proposition  \ref{prp:code}.\ref{prp:code9}
when $\mbox{\boldmath$A$}$ is the empty list in the right hand side,
and $o(P^{\prime})\prec o(P)$.
From IH we see that $\Gam_{0}\Rarw A_{0}$ is true.

The case when $A$ vanishes at a $(cut)$ is similar.
\\

\noindent
{\bf Case 2}.
The case when there exists a lower sequent of an $(ind)$ in the end-piece of $P$.
Let $P$ be the following:
\[
\infer*{\Gam_{0}\Rarw A_{0}}
{
 \infer[(ind)]
 {s:\Gam\Rarw C}
 {
  \infer*[Q_{0}(x)]{s_{0}: \Gam,\fal y\in x\, A(y)\Rarw A(x)}{}
 &
 \infer*[Q_{1}]{s_{1}: A(\bar{a}),\Gam\Rarw C}{}
 &
 s_{2}: \Gam\Rarw \bar{a}\in\bar{b}
 }
}
\]
If the formula $\bar{a}\in\bar{b}$ is false, i,e., $a\not\in b$, then replace $\bar{a}\in\bar{b}$ by $C$ in the succedents of the proof of 
$s_{2}: \Gam\Rarw \bar{a}\in\bar{b}$:
\[
\infer*{\Gam_{0}\Rarw A_{0}}
{
 \infer[(Rep)]
 {s:\Gam\Rarw C}
 {
 s_{2}: \Gam\Rarw C
 }
}
\]
We have $o(s;P)=(\gam_{0}\#\bar{1}\cdot 6)\times\bar{b}\# \gam_{1}\#\gam_{2}$
for $\gam_{i}=o(s_{i};P)$.
Since $o(s;P^{\prime})=o(s_{2};P^{\prime})=o(s_{2};P)=\gam_{2}\prec o(s;P)$,
we obtain $o(P^{\prime})\prec o(P)$.

Assume $\bar{a}\in\bar{b}$ is true, and let $P^{\prime}$ be the following:
{\small
\[
\infer*{\Gam_{0}\Rarw A_{0}}
{
 \infer{s: \Gam\Rarw C}
 {
  \infer{\Gam\Rarw A(\bar{a})}
  {
   \infer[(R\to, R\fal)]{\Gam\Rarw\fal y\in \bar{a}\, A(y)}
   {
    \infer[(ind)]
    {s^{\prime}:\Gam,z\in \bar{a}\Rarw A(z)}
    {
     \infer*[Q_{0}(x)]{s_{0}:\Gam,\fal y\in x\, A(y)\Rarw A(x)}{}
    &
   A(z),\Gam\Rarw A(z)
    &
    \Gam,z\in a\Rarw z\in a
    }
   }
   &
   \hspace{-19mm}
   \infer*[Q_{0}(a)]{s_{a}:  \Gam,\fal y\in \bar{a}\, A(y)\Rarw A(\bar{a})}{}
  }
  &
    \hspace{-10mm}
  \infer*[Q_{1}]{s_{1}:A(\bar{a}),\Gam\Rarw C}{}
 }
}
\]
}
where the proof $Q_{0}(a)$ is obtained from the subproof $Q_{0}(x)$ of $P$ by substituting the constant $\bar{a}$ 
for the eigenvariable $x$,
and renaming free variables for the pure variable condition for $P^{\prime}$.
The last two inference rules leading to ${s: \Gam\Rarw C}$ are $(cut)$'s.

It is easy to see that $\gam_{0}^{\prime}=o(s_{a};P^{\prime})\preceq o(s_{0};P)=\gam_{0}$
from $\bar{a}\prec mj(x)=\bar{V}$ for $a\in V$ and Proposition  \ref{prp:code}.\ref{prp:code8}.
We have $o(s^{\prime};P^{\prime})=(\gam_{0}\#\bar{1}\cdot 6)\times \bar{a}\#\bar{1}\cdot 4$.
By Proposition  \ref{prp:code}.\ref{prp:code5} we have
$(\gam_{0}\#\bar{1}\cdot 6)\times\bar{a}\#\bar{1}\cdot 6\#\gam_{0}^{\prime}\prec(\gam_{0}\# \bar{1}\cdot 6)\times\bar{b}$.
Hence 
we obtain 
$o(s;P^{\prime})=(\gam_{0}\#\bar{1}\cdot 6)\times\bar{a}\#\bar{1}\cdot 6\#\gam_{0}^{\prime}\#\gam_{1}\prec
(\gam_{0}\#\bar{1}\cdot 6)\times\bar{b}\# \gam_{1}\#\gam_{2}=o(s;P)$.
This yields $o(P^{\prime})\prec o(P)$.
\\

\noindent
In the following two cases inference rules introducing $Q$-free formulae and $(cut)$ with $Q$-free cut formulae
are pushed down to the end of proofs.
\\

\noindent
{\bf Case 3}.
The case when there exists a lower sequent of an explicit inference rule in the end-piece of $P$,
where an inference rule $J$ is explicit in $P$ iff its major (principal) formula is in the antecedents (succedents) of any sequent below it
when the formula is in the antecedent (succedent) of the lower sequent of $J$, resp.

Let $J$ be such an inference rule. $J$ is one of the inference rules 
$(L\lor)$, $(R\lor)$, $(L\land)$, $(R\land)$, $(L\to)$, $(R\to)$, 
$(L\exists)$, $(R\exists)$, $(L\forall)$, and $(R\forall)$, but neither of $(LQ)$ and $(RQ)$, 
since the fixed point predicate $Q$ does not occur in the end sequent of $P$.

Consider the cases when $J$ is either an $(R\fal)$ or an $(L\to)$.
For the first case let $P$ be the following:
\[
\infer*{\Gam_{0}\Rarw \fal x\, A(x)}
{
 \infer[(R\fal)]
 {s:\Gam\Rarw \fal x\,A(x)}
 {
 \infer*[Q(y)]{s_{0}:\Gam\Rarw A(y)}{}
 }
}
\]
For each $a\in V$, let $P_{a}$ be the following:
\[
\infer*{\Gam_{0}\Rarw A(\bar{a})}
{
 \infer[(Rep)]
 {s:\Gam\Rarw A(\bar{a})}
 {
 \infer*[Q(a)]{s_{a}:\Gam\Rarw A(\bar{a})}{}
 }
}
\]
Since $o(s;P_{a})=o(s_{a};P_{a})\preceq o(s_{0};P)\prec o(s;P)$,
we have $o(P_{a})\prec o(P)$.
By IH $\Gam_{0}\Rarw A(\bar{a})$ is true for any $a\in V$.
Hence so is $\Gam_{0}\Rarw\fal x\, A(x)$.

For the second case let $P$ be the following:
\[
\infer*[Q]{\Gam_{0}\Rarw A_{0}}
{
 \infer[(L\to)]
 {s:\Gam, B\to C\Rarw A_{1}}
 {
 \Gam,B\to C\Rarw B
 &
 \Gam,B\to C,C\Rarw A_{1}
 }
}
\]
where $(B\to C)\in\Gam_{0}$.

Let $P_{C}$ be the following:
\[
\infer*{\Gam_{0},C \Rarw A_{0}}
{
 \infer[(Rep)]
 {s:\Gam, B\to C,C\Rarw A_{1}}
 {
 \Gam,B\to C,C\Rarw A_{1}
 }
}
\]
Since $o(s;P^{\prime})\prec o(s;P)$, we obtain $o(P_{C})\prec o(P)$, and 
$\Gam_{0},C \Rarw A_{0}$ is true by IH.

Next let $P_{B}$ be the following:
\[
\infer*[Q_{B}]{\Gam_{0} \Rarw B}
{
 \infer[(Rep)]
 {s:\Gam, B\to C\Rarw B}
 {
 \Gam,B\to C\Rarw B
 }
}
\]
where the trunk $Q_{B}$ is obtained from the trunk $Q$ of $P$ as follows.
If in $Q$, $A_{1}$ vanishes as a cut formula,
\[
\infer[(chain)]{\mbox{\boldmath$\Gam$}_{1},\Gam_{1},B\to C,\Del\Rarw D}
{
\mbox{\boldmath$\Gam$}_{1}\Rarw \mbox{\boldmath$A$}
&
\infer*{\Gam_{1}, B\to C\Rarw A_{1}}{s:\Gam, B\to C\Rarw A_{1}}
&
\Del,\mbox{\boldmath$A$},A_{1}\Rarw D
}
\]
then this part turns to
\[
\infer[(Rep)]{\mbox{\boldmath$\Gam$}_{1}, \Gam_{1},B\to C,\Del\Rarw B}
{
\infer*{\Gam_{1}, B\to C\Rarw B}{s:\Gam, B\to C\Rarw B}
}
\]
This pruning step is iterated when $D$ vanishes below.
Clearly we have $o(P_{B})\prec o(P)$, and $\Gam_{0} \Rarw B$ is true by IH.

Since both $\Gam_{0},C \Rarw A_{0}$ and $\Gam_{0} \Rarw B$ are true, and $(B\to C)\in\Gam_{0}$,
so is $\Gam_{0}\Rarw A_{0}$.
\\

\noindent
{\bf Case 4}.
The case when there exists a cut formula $A_{1}$ in the end-piece of $P$ such that $A_{1}$ is a
$Q$-free formula.

Let $P$ be the following:
\[
\infer*{\Gam_{0}\Rarw A_{0}}
{
 \infer[(chain)]{\mbox{\boldmath$\Gam$},\Gam_{1},\Del\Rarw C}
 {
  \mbox{\boldmath$\Gam$}\Rarw \mbox{\boldmath$A$}
  &
  \Gam_{1}\Rarw A_{1}
  &
  \Del,\mbox{\boldmath$A$},A_{1}\Rarw C
  }
}
\]
Let $P_{r}$ be the following which is obtained from $P$ as for $P_{B}$ in the {\bf Case 3}.
\[
\infer*{\Gam_{0}\Rarw A_{1}}
{
 \infer[(Rep)]{\mbox{\boldmath$\Gam$},\Gam_{1},\Del\Rarw A_{1}}
 {
  \Gam_{1}\Rarw A_{1}
  }
}
\]
And let $P_{\ell}$ be the following:
\[
\infer*{\Gam_{0},A_{1}\Rarw A_{0}}
{
 \infer[J]{\mbox{\boldmath$\Gam$},\Gam_{1},\Del,A_{1}\Rarw C}
 {
  \mbox{\boldmath$\Gam$}\Rarw \mbox{\boldmath$A$}
  &
  \Del,\mbox{\boldmath$A$},A_{1}\Rarw C
  }
}
\]
where $J$ denotes two consecutive $(E)$'s if $\mbox{\boldmath$A$}$ is the empty list, and a $(chain)$ otherwise.

Obviously both $P_{r}$ and $P_{\ell}$ are restricted, and $o(P_{r}),o(P_{\ell})\prec o(P)$.
IH says that both $\Gam_{0}\Rarw A_{1}$ and $\Gam_{0},A_{1}\Rarw A_{0}$ are true.
Hence so is $\Gam_{0}\Rarw A_{0}$.

The case when $A_{1}$ is a cut formula of a $(cut)$ is similar.
\\

\noindent
{\bf Case 5}.
The case when there exists a $(cut)\, J_{0}$ in the end-piece of $P$ such that for its lower sequent $s: \Gam,\Del\Rarw C$
and cut formula $A$,
$h(s;P)>d:=dp(A)>0$.
Let $J$ be the uppermost $(E)$ below $J_{0}$.
Note that here is no $(chain)$ between $J_{0}$ and $J$ since $P$ is height-normal.
Let $P$ be the following.
\[
\infer*{\Gam_{0}\Rarw A_{0}}
{
 \infer[(E)\, J]{u: \Gam_{1}\Rarw C_{1}}
 {
  \infer*{t: \Gam_{1}\Rarw C_{1}}
  {
   \infer[(cut)\, J_{0}]{ s: \Gam,\Del\Rarw C}
   {
    s_{1}:\Gam\Rarw A
    &
    s_{2}:\Del,A\Rarw C
    }
   }
  }
}
\]

Let $P^{\prime}$ be obtained from $P$ by lowering the $(cut)\, J_{0}$ below the $(E)\, J$:

\[
\infer*{\Gam_{0}\Rarw A_{0}}
{
  \infer[(cut)]{u: \Gam_{1}\Rarw C_{1}}
  {
    \infer[(E)]{u_{1}: \Gam_{1}\Rarw A}
    {
      \infer*{t_{1}: \Gam_{1}\Rarw A}
      {
       \infer[(Rep)]{s: \Gam,\Del\Rarw A}
        {
          s_{1}:\Gam\Rarw A
        }
      }
     }
  &
    \infer[(E)]{u_{2}: \Gam_{1}, A\Rarw C_{1}}
    {
      \infer*{t_{2}: \Gam_{1}, A\Rarw C_{1}}
      {
        \infer[(Rep)]{s:\Gam,\Del, A \Rarw C}
        {
        s_{2}:\Del,A \Rarw C
        }
       }
     }
   }
}
\]

Let $\alp_{i}=o(s_{i};P)=o(s_{i};P^{\prime})$ for $i=1,2$.
Then for some $\bet$,
$o(t;P)=\bet\#\alp_{1}\#\alp_{2}$, and $o(u;P)=\Ome^{\bet\#\alp_{1}\#\alp_{2}}$.
On the other side for some $\bet^{\prime}\preceq\bet$, 
$o(t_{1};P^{\prime})=\bet^{\prime}\#\alp_{1}$.
The case $\bet^{\prime}\prec\bet$ happens when a pruning is performed. 
Also $o(t_{2};P^{\prime})=\bet\#\alp_{2}$.
Hence $o(u_{1};P^{\prime})\preceq \Ome^{\bet\#\alp_{1}}$ and $o(u_{2};P^{\prime})=\Ome^{\bet\#\alp_{2}}$.
Now we claim that $o(u;P^{\prime})\preceq \Ome^{\bet\#\alp_{1}}\#\Ome^{\bet\#\alp_{2}}\prec \Ome^{\bet\#\alp_{1}\#\alp_{2}}=o(u;P)$, which follows from Proposition  \ref{prp:code}.\ref{prp:code10}.

Hence $o(P^{\prime})\prec o(P)$, and we see that $\Gam_{0}\Rarw A_{0}$ is true from IH.
\\

\noindent
In the following cases, adjacent $(cut)$'s are first collected into $(chain)$, {\bf Case 6}.
This as well as the analysis of strictly positive cut formula in {\bf Case 9} prolongs $(chain)$.
In {\bf Case 7}, $(cut)$ with strictly positive cut formula is replaced by $(chain)$, thereby $(chain)$ is introduced in proofs.
\\

\noindent
{\bf Case 6}.
The case when there exists a $(cut)\, J_{0}$ in the end-piece of $P$ such that its lower sequent $s: \Gam_{1},\Del_{1}\Rarw C$
is the rightmost upper sequent of a $(chain)\, J$.
Let $P$ be the following with $\Del=\Del_{0}\cup\Del_{1}$:
\[
\infer*{\Gam_{0}\Rarw A_{0}}
{
 \infer[(chain)\, J]{s: \mbox{\boldmath$\Gam$},\Del\Rarw C}
 {
  \mbox{\boldmath$s$}: \mbox{\boldmath$\Gam$}\Rarw \mbox{\boldmath$A$}
  &
  \infer[(cut)\, J_{0}]{\Del,\mbox{\boldmath$A$}\Rarw C}
  {
   s_{0}: \Del_{0},\mbox{\boldmath$A$}\Rarw A_{0}
   &
  t: \Del_{1},\mbox{\boldmath$A$},A_{0}\Rarw C
   }
  }
}
\]
Since $P$ is height-normal, we have $2=h(t;P)\geq dp(A_{0})$.
On the other side $h(t;P)\leq dp(A_{0})$ by virtue of {\bf Case 5}.
Hence $dp(A_{0})=2$, i.e., the predicate $Q$ occurs in $A_{0}$ and $A_{0}$ is strictly positive.

Let $P^{\prime}$ be the following:
\[
\infer*{\Gam_{0}\Rarw A_{0}}
{
 \infer[(chain)]{s: \mbox{\boldmath$\Gam$},\Del\Rarw C}
 {
  \mbox{\boldmath$s$}: \mbox{\boldmath$\Gam$}\Rarw \mbox{\boldmath$A$}
  &
  \infer[(chain)]{\mbox{\boldmath$\Gam$},\Del_{0}\Rarw A_{0}}
   {
    \mbox{\boldmath$s$}: \mbox{\boldmath$\Gam$}\Rarw \mbox{\boldmath$A$}
   &
     s_{0}: \Del_{0},\mbox{\boldmath$A$}\Rarw A_{0}
    }
  &
  \infer[(Rep)]{\Del_{1},\mbox{\boldmath$A$},A_{0}\Rarw C}
  {
  t: \Del_{1},\mbox{\boldmath$A$},A_{0}\Rarw C
   }
  }
}
\]
Observe that $2=h(s_{0};P^{\prime})=h(s_{0};P)=h(t;P)=h(t;P^{\prime})$.
Let $\mbox{\boldmath$\alp$}=o(\mbox{\boldmath$s$};P)=o(\mbox{\boldmath$s$};P^{\prime})$,
$\alp_{0}=o(s_{0};P)=o(s_{0};P^{\prime})$ and $\bet=o(t;P)=o(t;P^{\prime})$.
Then $o(s;P)=\Ome_{2}(\alp_{0}\#\bet)\sum\mbox{\boldmath$\alp$}$ and
$o(s;P^{\prime})=\Ome_{2}(\bet)(\sum\mbox{\boldmath$\alp$}\#(\Ome_{2}(\alp_{0})\sum\mbox{\boldmath$\alp$}))$.
$o(s;P^{\prime})\prec o(s;P)$ is seen from Proposition \ref{prp:code}.\ref{prp:code10}.
Hence $o(P^{\prime})\prec o(P)$, and we see that $\Gam_{0}\Rarw A_{0}$ is true from IH.
\\

\noindent
{\bf Case 7}.
The case when there exists a $(cut)$ with a strictly positive cut formula $A$ in the end-piece of $P$.
Let $J$ be a lowest such $(cut)$.
By virtue of {\bf Case 5} we have $h(t;P)=dp(A)=2$, and by {\bf Case 6} there is no rightmost upper sequent of any $(chain)$
below $J$.
Hence there are two consecutive $(E)$'s below $J$ by {\bf Case 4} and {\bf Case 5}.
Furthermore the two consecutive $(E)$'s is immediately below the lowest $J$, i.e.,
there is no left upper sequent of any $(chain)$ between $J$ and $(E)$'s since $P$ is height-normal.
Let $P$ be the following:
\[
\infer*{\Gam_{0}\Rarw A_{0}}
{
 \infer[(E)^{2}]{s:\Gam,\Del\Rarw C}
 {
  \infer[(cut)\, J]{t: \Gam,\Del\Rarw C}
   {
    \infer*{u_{0}: \Gam\Rarw A}{}
   &
    \infer*{u_{1}: \Del,A\Rarw C}{}
    }
 }
}
\]
Let $P^{\prime}$ be the following:
\[
\infer*{\Gam_{0}\Rarw A_{0}}
{
\infer[(chain)]{s^{\prime}:\Gam,\Del\Rarw C}
{
  \infer[(E)^{2}]{s_{0}:\Gam,\Del\Rarw A}
  {
   \infer[(Rep)]{t_{0}: \Gam,\Del\Rarw A}
   {
    \infer*{u_{0}^{\prime}: \Gam\Rarw A}{}
    }
  }
&
  \infer[(Rep)]{t_{1}:\Gam,\Del,A\Rarw C}
  {
   \infer*{u_{1}^{\prime}: \Del,A\Rarw C}{}
  }
 }
}
\]
We have $h(s;P)=h(s^{\prime};P^{\prime})=h(s_{0};P^{\prime})=0$ and
$h(t_{0};P^{\prime})=h(t_{1};P^{\prime})=h(t;P)=2$.
Let $\alp_{i}=o(u_{i};P)$ for $i=0,1$.
Then
$o(u_{0}^{\prime};P^{\prime})=\alp_{0}$, $o(t_{1};P^{\prime})=o(u_{1}^{\prime};P^{\prime})=\alp_{1}$, and
$o(s_{0};P^{\prime})=\Ome_{2}(\alp_{0})$.
Hence $o(s^{\prime};P^{\prime})=\Ome_{2}(\alp_{1}) \Ome_{2}(\alp_{0})=\Ome^{\Ome^{\alp_{1}}\#\Ome^{\alp_{0}}}
\prec \Ome_{2}(\alp_{0}\#\alp_{1})=o(s;P)$ by Proposition \ref{prp:code}.\ref{prp:code10}.
Therefore $o(P^{\prime})\prec o(P)$, and by IH $\Gam_{0}\Rarw A_{0}$ is true.
\\

\noindent
By virtue of {\bf Case 1}-{\bf Case 3} we can assume that any topmost sequent in the end-piece of $P$ is a lower sequent
of an implicit inference rule other than $(ind)$, $(cut)$, $(chain)$, $(Rep)$ and $(E)$
such that the fixed point predicate $Q$ occurs in its major formula.
Call temporarily such an inference rule {\it boundary\/} of $P$ if its lower sequent is in the end-piece, but not its upper sequents.
We then claim that there is an inference $J$ such that $J$ is either a $(cut)$ or a $(chain)$,
and one of its cut formula $A$ comes from major formulae of boundaries.
\[
\infer*{\Gam_{0}\Rarw A_{0}}
{
 \infer[J]{\mbox{\boldmath$\Gam$}, \Gam,\Del\Rarw C}
 {
 \mbox{\boldmath$\Gam$}\Rarw \mbox{\boldmath$A$}
 &
 \infer*{\Gam\Rarw A}
  {
  \infer[J_{\ell}]{\Gam_{1}\Rarw A}{}
   }
 &
 \infer*{\Del,\mbox{\boldmath$A$},A\Rarw C}
  {
   \infer[J_{r}]{\Del_{1},A\Rarw C_{1}}{}
   }
 }
}
\]
where both $J_{\ell}$ and $J_{r}$ are boundaries, $A$ in their lower sequents are their major formulae, and
the formula $A$ is in the succednets [antecedents] of any sequents between $J_{\ell}$ and $J$
[between $J_{r}$ and $J$], resp.

The claim is seen as in \cite{PT2} (the existence of a suitable cut).

In what follows pick such rules $J$, $J_{\ell}$ and $J_{r}$ with the formula $A$, which is a cut formula of $J$.
By virtue of {\bf Case 7}, $J$ is a $(cut)$ iff $dp(A)>2$.
\\

\noindent
{\bf Case 8}.
The case when $dp(A)>2$ and $J$ is a $(cut)$.
For example consider the case when $A$ is a formula $\fal x\, D(x)$.
Let $P$ be the following:
\[
\infer*{\Gam_{0}\Rarw A_{0}}
{
\infer[(E)\, J_{0}]{t: \Gam_{2}\Rarw B}
 {
  \infer*{v: \Gam_{2}\Rarw B}
  {
   \infer[(cut)\, J]{s: \Gam,\Del\Rarw C}
   {
    \infer*{u_{0}:\Gam\Rarw \fal x\, D(x)}
    {
      \infer[(R\fal)\, J_{\ell}]{\Gam_{1}\Rarw \fal x\, D(x)}{\Gam_{1}\Rarw D(y)}
     }
   &
     \infer*{u_{1}:\Del,\fal x\, D(x)\Rarw C}
     {
      \infer[(L\fal)\, J_{r}]{\Del_{1},\fal x\, D(x)\Rarw C_{1}}{\Del_{1},\fal x\, D(x),D(a)\Rarw C_{1} }
     }
    }
   }
 }
}
\]
By virtue of {\bf Case 5} we can assume that $h(s;P)=dp(\fal x\, D(x))=d+1$ with $d=dp(D(a))>2$.
$J_{0}$ denotes the uppermost $(E)$ below $J$ with $h(t;P)=d$.

Let $P^{\prime}$ be the following:
\[
\infer*{\Gam_{0}\Rarw A_{0}}
{
 \infer[(cut)]{t^{\prime}:\Gam_{2}\Rarw B}
 {
  \infer[(E)]{t_{\ell}: \Gam_{2}\Rarw D(a)}
  {
   \infer*{v_{\ell}: \Gam_{2}\Rarw D(a)}
   {
     \infer[(Rep)]{s_{\ell}: \Gam,\Del\Rarw D(a)}
    {
     \infer*{u_{0}^{\prime}:\Gam\Rarw D(a)}
      {
       \infer[(Rep)]{\Gam_{1}\Rarw D(a)}{\Gam_{1}\Rarw D(a)}
       }
     }
    }
  }
&
 \infer[(E)]{t_{r}: \Gam_{2},D(a)\Rarw B}
 {
  \infer*{v_{r}: \Gam_{2},D(a)\Rarw B}
  {
   \infer[(cut)]{s_{r}: \Gam,\Del,D(a)\Rarw C}
   {
    \infer*{u_{0}:\Gam\Rarw \fal x\, D(x)}{}
   &
       \infer*{u_{1}^{\prime}: \Del,\fal x\, D(x),D(a)\Rarw C}
       {
        \infer[(Rep)]{\Del_{1},\fal x\, D(x),D(a)\Rarw C_{1}}{\Del_{1},\fal x\, D(x),D(a)\Rarw C_{1} }
       }
     }
    }
  }
}
}
\]
We have $o(s;P)=\alp_{0}\#\alp_{1}$ where $\alp_{i}=o(u_{i};P)$ for $i=0,1$.
On the other hand we have $o(s_{\ell};P^{\prime})=\alp_{0}^{\prime}=o(u_{0}^{\prime};P^{\prime})\prec\alp_{0}=o(u_{0};P^{\prime})\prec o(s_{r};P^{\prime})$
and $\alp_{1}^{\prime}=o(u_{1}^{\prime};P^{\prime})\prec \alp_{1}$.
Hence $o(s_{\ell};P^{\prime})\prec o(s_{r};P^{\prime})\prec o(s;P)$, and
$o(u_{\ell};P^{\prime})\prec o(u_{r};P^{\prime})\prec o(u;P)$.
Thus for $o(t_{\ell};P^{\prime})=\Ome^{o(u_{\ell};P^{\prime})}, o(t_{r};P^{\prime})=\Ome^{o(u_{r};P^{\prime})}$, and $\Ome^{o(u;P)}=o(t;P)$,
we obtain $o(t^{\prime};P^{\prime})=o(t_{\ell};P^{\prime})\#o(t_{r};P^{\prime})\prec o(t;P)$.
Therefore $o(P^{\prime})\prec o(P)$, and by IH $\Gam_{0}\Rarw A_{0}$ is true.

The other cases are seen similarly.
\\

\noindent
{\bf Case 9}.
The case when $dp(A)=2$ and $J$ is a $(chain)$.

First consider the case when $A$ is an implicational formula $D\to E$,
where $E$ is strictly positive and $D$ is $Q$-free.
Let $P$ be the following:
{\small
\[
\infer*{\Gam_{0}\Rarw A_{0}}
{
 \infer[(chain)\, J]{s_{0}: \mbox{\boldmath$\Gam$},\Gam,\Del\Rarw C}
 {
 \mbox{\boldmath$s$}: \mbox{\boldmath$\Gam$}\Rarw \mbox{\boldmath$A$}
 &
 \infer*{s_{1}: \Gam\Rarw D\to E}
  {
  \infer[(R\to)\, J_{\ell}]{s_{3}: \Gam_{1}\Rarw D\to E}{s_{4}: \Gam_{1},D\Rarw E}
   }
 &
 \infer*{s_{2}: \Del,\mbox{\boldmath$A$},D\to E\Rarw C}
  {
   \infer[(L\to)\, J_{r}]{s_{5}: \Del_{1},D\to E\Rarw C_{1}}
   {
   s_{6}: \Del_{1}, D\to E\Rarw D
   &
   s_{7}:  \Del_{1},D\to E, E\Rarw C_{1}
    }
   }
 }
}
\]
}
Let $P_{\ell}$ be the following:
\[
\infer*{\Gam_{0}\Rarw D}
{
  \infer[(chain)]{s_{\ell}: \mbox{\boldmath$\Gam$},\Gam,\Del\Rarw D}
  {
   \mbox{\boldmath$\Gam$}\Rarw \mbox{\boldmath$A$}
   &
    \infer*{\Gam\Rarw D\to E}{}
   &
    \infer*{s_{2\ell}:\Del,\mbox{\boldmath$A$},D\to E\Rarw D}
    {
     \infer[(Rep)]{\Del_{1},D\to E\Rarw D}
     {
     s_{6}:\Del_{1}, D\to E\Rarw D
      }
    }
   }
 }
\]
Let $P_{r}$ be the following:
\[
\infer*{\Gam_{0},D\Rarw A_{0}}
{
  \infer[(chain)]{s_{r}: \mbox{\boldmath$\Gam$},\Gam,\Del,D\Rarw C}
  {
  \mbox{\boldmath$\Gam$}\Rarw \mbox{\boldmath$A$}
    &
    \infer*{\Gam\Rarw D\to E}{}
   &
    \infer*{s_{1r}:\Gam,D\Rarw  E}
    {
    \infer[(Rep)]{\Gam_{1},D\Rarw  E}{s_{4}:\Gam_{1},D\Rarw E}
     }
  &
    \infer*{s_{2r}:\Del,\mbox{\boldmath$A$},D\to E,E\Rarw C}
    {
     \infer[(Rep)]{\Del_{1},D\to E,E\Rarw C_{1}}
     {
     s_{7}:\Del_{1},D\to E, E\Rarw C_{1}
     }
   }
   }
 }
\]
Let $\alp_{i}=o(s_{i};P)$ for $i=4,6,7$.
In $P$, $o(s_{3};P)=\alp_{4}\# \bar{1}$, $o(s_{5};P)=\alp_{6}\#\alp_{7}$, and
$o(s_{0};P)=\Ome_{2}(o(s_{2};P))(\sum\mbox{\boldmath$\alp$}\#o(s_{1};P))$
for $ \mbox{\boldmath$\alp$}= o(\mbox{\boldmath$s$};P)$.
On the other side in $P_{\ell}$ and $P_{r}$,
$\alp_{6}=o(s_{6};P_{\ell})$, $\alp_{4}=o(s_{4};P_{r})$ and $\alp_{7}=o(s_{7};P_{r})$, and hence 
$o(s_{2\ell};P_{\ell})\prec o(s_{2};P)$,
$o(s_{1r};P_{r})\prec o(s_{1};P)$ and 
$o(s_{2r};P_{r})\prec o(s_{2};P)$.
Moreover
\\
$o(s_{\ell};P^{\prime})=\Ome_{2}(o(s_{2\ell};P_{\ell}))(\sum\mbox{\boldmath$\alp$}\#o(s_{1};P))$
and
\\
$o(s_{r};P_{r})=\Ome_{2}(o(s_{2r};P_{r}))(\sum\mbox{\boldmath$\alp$}\#o(s_{1};P)\#o(s_{1r};P_{r}))$.

We see $o(s_{\ell};P_{\ell}), o(s_{r};P_{r})\prec o(s_{0};P)$ from Proposition \ref{prp:code}.\ref{prp:code11}.
From these we see that $o(P_{\ell}),o(P_{r})\prec o(P)$, and by IH 
both $\Gam_{0}\Rarw D$ and $\Gam_{0},D\Rarw A_{0}$ are true.
Therefore $\Gam_{0}\Rarw A_{0}$ is true.

Next consider the case when $A\equiv Q(a)$ for the fixed point predicate $Q$.
Let $P$ be the following:
\[
\infer*{\Gam_{0}\Rarw A_{0}}
  {
  \infer[(chain)\, J]{s_{0}: \mbox{\boldmath$\Gam$},\Gam,\Del\Rarw C}
   {
    \mbox{\boldmath$s$}: \mbox{\boldmath$\Gam$}\Rarw \mbox{\boldmath$A$}
  &
    \infer*{s_{1}: \Gam\Rarw Q(a)}
     {
     \infer[(RQ)\, J_{\ell}]{\Gam_{1}\Rarw Q(a)}{\Gam_{1}\Rarw \mathcal{Q}(Q,a)}
      }
   &
    \infer*{s_{2}: \Del,\mbox{\boldmath$A$}, Q(a)\Rarw C}
    {
     \infer[(LQ)\, J_{r}]{\Del_{1},Q(a)\Rarw C_{1}}{\Del_{1},Q(a),\mathcal{Q}(Q,a)\Rarw C_{1} }
     }
   }
}
\]

Let $P^{\prime}$ be the following:
\[
\infer*{\Gam_{0}\Rarw A_{0}}
  {
  \infer[(chain)]{s_{0}: \mbox{\boldmath$\Gam$},\Gam,\Del\Rarw C}
   {
    \mbox{\boldmath$s$}: \mbox{\boldmath$\Gam$}\Rarw \mbox{\boldmath$A$}
  &
     \infer*{s_{1}: \Gam\Rarw Q(a)}{}
   &
    \infer*{s_{1}^{\prime}: \Gam\Rarw \mathcal{Q}(Q,a)}
     {
     \infer[(Rep)]{\Gam_{1}\Rarw \mathcal{Q}(Q,a)}{\Gam_{1}\Rarw \mathcal{Q}(Q,a)}
      }
   &
    \infer*{s_{2}: \Del,\mbox{\boldmath$A$}, Q(a),\mathcal{Q}(Q,a)\Rarw C}
    {
     \infer[(Rep)]{\Del_{1},Q(a),\mathcal{Q}(Q,a)\Rarw C_{1}}{\Del_{1},Q(a),\mathcal{Q}(Q,a)\Rarw C_{1} }
     }
   }
}
\]
We have
$o(s_{0};P)=\Ome_{2}(o(s_{2};P))(\sum\mbox{\boldmath$\alp$}\#o(s_{1};P))$
for $ \mbox{\boldmath$\alp$}= o(\mbox{\boldmath$s$};P)$.
Also 
\\
$o(\mbox{\boldmath$s$};P^{\prime})=o(\mbox{\boldmath$s$};P)$,
$o(s_{1}^{\prime};P^{\prime})\prec o(s_{1};P)=o(s_{1};P^{\prime})$, 
$o(s_{2};P^{\prime})\prec o(s_{2};P)$, and
$o(s_{0};P^{\prime})=\Ome_{2}(o(s_{2};P^{\prime}))(\sum\mbox{\boldmath$\alp$}\#o(s_{1};P)\#o(s_{1}^{\prime};P^{\prime}))$.
Hence
$o(s_{0};P^{\prime})\prec o(s_{0};P)$ from Proposition \ref{prp:code}.\ref{prp:code11}.
Therefore $o(P^{\prime})\prec o(P)$, and by IH $\Gam_{0}\Rarw A_{0}$ is true.

The other cases are seen similarly.
This completes a proof of Lemma \ref{lem:main}, and of Theorem \ref{th:consvintfix}.

\end{document}